\documentclass[oneside,final]{csri24}

\usepackage{xcolor}
\usepackage{xurl}
\usepackage[normalem]{ulem}

\usepackage[top=1.355in,
			bottom=1.355in,
			left=1.5in,
			right=1.5in]{geometry}
\usepackage{amsfonts,
            amsmath,
            graphicx,
            url}
\usepackage{mathtools}
\usepackage{booktabs}
\usepackage{caption}
\usepackage{subcaption}
\usepackage{boldtensors}
\usepackage{makecell}
\usepackage{placeins}

\title{Transmission Conditions for the Non-Overlapping Schwarz Coupling of Full Order and Operator Inference Models}

\author{Cameron Rodriguez\thanks{Columbia University,  cjr2210@columbia.edu} \and Irina Tezaur\thanks{Sandia National Laboratories, ikalash@sandia.gov} \and Alejandro Mota\thanks{Sandia National Laboratories, amota@sandia.gov} \and Anthony Gruber\thanks{Sandia National Laboratories, adgrube@sandia.gov}  \and Eric Parish\thanks{Sandia National Laboratories, ejparis@sandia.gov}  \and Christopher Wentland\thanks{Sandia National Laboratories, crwentl@sandia.gov} }

\begin{document}

\maketitle


\begin{abstract}

This work investigates transmission conditions for the domain decomposition-based coupling of subdomain-local models using the non-overlapping Schwarz alternating method (NO-SAM). Building on prior efforts involving overlapping SAM (O-SAM), we formulate and assess two NO-SAM variants, based on alternating Dirichlet–Neumann and Robin–Robin transmission conditions.  For the subdomain-local models, we consider a mix of full order models (FOMs) and non-intrusive reduced order models (ROMs) constructed via an emerging model reduction technique known as operator inference (OpInf).  Of particular novelty is the first application of NO-SAM to couple non-intrusive OpInf ROMs with each other, and with FOMs. Numerical studies on a one-dimensional linear elastic wave propagation benchmark problem demonstrate that transmission condition choice and parameter tuning significantly impact convergence rate, accuracy, and stability. Robin–Robin coupling often yields faster convergence than alternating Dirichlet–Neumann, though improper parameter selection can induce spurious oscillations at subdomain interfaces. For FOM–OpInf and OpInf–OpInf couplings, sufficient modal content in the ROM basis improves accuracy and mitigates instability, in some cases outperforming the coupled FOM–FOM reference solutions in both accuracy and efficiency.  
These findings highlight NO-SAM’s potential for enabling flexible, non-intrusive, and efficient multi-model coupling across independently meshed subdomains, while emphasizing the need for careful interface condition design  in higher-dimensional and predictive settings.

\end{abstract}

\section{Introduction} \label{sec:intro}


Despite advancements in algorithms and high-performance computing (HPC), high-fidelity modeling and simulation (mod/sim) of complex multi-scale and multi-physics systems remains a major challenge due to: (i) the time-intensive process of generating high-quality meshes, identified as ``the single biggest bottleneck in [mod/sim] analyses" \cite{SandiaLabNews}, and (ii) the substantial computational cost and long runtimes required for performing high-fidelity simulations of complex nonlinear phenomena.  While data-driven reduced order models (ROMs) have promised to mitigate the runtime burden of mod/sim, these models face their own set of challenges related to robustness, stability and accuracy, along with lengthy implementations and a lack of systematic refinement mechanisms. 

The present work advances a computational technique that promises to mitigate both hurdles described above: the Schwarz alternating method (SAM) \cite{Schwarz:1870}. The key idea behind SAM is to decompose the physical domain on which a given partial differential equation (PDE) is posed into smaller subdomains, and then to solve the governing PDE by iterating between smaller, subdomain-local problems, while exchanging boundary condition information to ensure compatibility across subdomain interfaces. 
Over the past decade, SAM has been shown to enable concurrent domain decomposition- (DD-)based coupling of subdomain-local full order models (FOMs) with each other \cite{Mota:2017, Mota:2022, Mota:2025}, as well as with a variety of data-driven surrogate models (e.g., intrusive projection-based ROMs \cite{Wentland:2025, Barnett:2022Schwarz}, physics-informed neural networks (PINNs) \cite{Snyder:2023}, non-intrusive operator inference (OpInf) models \cite{Moore:2024}). In these papers, SAM exhibits a number of advantages over competing multi-scale coupling methods. Most notably: it is concurrent and minimally-intrusive, can couple non-conformal meshes with different element topologies, and allows each subdomain to be advanced with different time integrators and time steps for dynamic problems. Importantly, SAM achieves this without introducing non-physical artifacts into the solution. 

The majority of our past research has focused on the overlapping version of SAM (O-SAM), which is known to converge to the monolithic solution of a well-posed PDE when simple Dirichlet transmission boundary conditions (BCs) are prescribed on the subdomain interfaces, and the overlap region is non-empty \cite{Schwarz:1870, Mota:2017, Mota:2022}.  Herein, we explore the performance of two variants of non-overlapping SAM (NO-SAM): alternating Dirichlet–Neumann NO-SAM (DN-NO-SAM) \cite{Zanolli:1987} and Robin–Robin NO-SAM (RR-NO-SAM) \cite{Lions:1990}.  Our interest in the NO version of SAM stems from the fact that this SAM variant provides maximal meshing flexibility, as subdomains can be created and meshed completely independently of each other without having to worry about the size/shape of the overlap region.  These independently-meshed subdomains can then be added to a library of archetypes, which can be glued together in arbitrary combinations using NO-SAM when creating new systems/geometries. In addition to formulating and assessing NO-SAM for FOM-FOM coupling, we also present the first (to our knowledge) formulation of NO-SAM for the DD-based coupling of non-intrusive OpInf ROMs \cite{willcox2016opinf} with each other and/or with FOMs. OpInf is a completely non-intrusive ROM method: since OpInf works by assuming a functional form for the ROM in terms of to-be-learned reduced operators and solving an optimization problem offline for these operators, it does not require access to the underlying FOM code like traditional (intrusive) projection-based ROM construction methods.  

The work described here is an addition to a growing body of literature on DD-based OpInf-OpInf and ROM-FOM couplings.  Our SAM approach is most similar to the recent work by Farcas et al. \cite{Farcas:2023}, which develops a method for performing DD-based coupling of subdomain-local OpInf ROMs by learning the appropriate reduced operators responsible for the coupling.  The present work differs from that of Farcas et al. in two important ways: (i) in \cite{Farcas:2023}, each subdomain problem is solved once rather than by performing an iteration to convergence as done within our Schwarz framework, which may be insufficient for certain classes of nonlinear problems \cite{SchwarzOpInfPaper:2025}, and (ii) the approach in \cite{Farcas:2023} is limited to overlapping domain decompositions (DDs). While other Schwarz-like methods for FOM-FOM and FOM-ROM couplings exist in the literature, the majority of these methods focus on the coupling of intrusive projection-based ROMs \cite{Corigliano:2015, Kerfriden:2013, Kerfriden:2012, Radermacher:2014} or neural network- (NN-)based ROMs \cite{Wang:2022, Wang:2025, LiD3M, LiDeepDDM, Snyder:2023}.  A variety of other FOM-FOM and FOM-ROM coupling methods, e.g., Lagrange multipliers, optimization-based coupling, coupling via boundary maps, each with their own advantages and disadvantages, have also been proposed over the years.  The interested reader is referred to \cite{  Cinquegrana:2011,  deCastro:2023, deCastro:2022,Iollo:2022,Maier:2014, Prusak:2023, Prusak:2024, Discacciati:2024, Hawkins:2024, Bochev:2024, Gkimisis:2025, Smetana:2023, Discacciati:2024, Hoang:2021, Chung:2024, Taddei:2024, Huang:2024, Ebrahimi:2024} and the references therein for more details. To the best of the authors' knowledge, only two of these references, namely \cite{Gkimisis:2025} and \cite{Farcas:2023}, consider the DD-based coupling of subdomain-local non-intrusive OpInf ROMs of the flavor considered herein \cite{willcox2016opinf}.

The remainder of this paper is organized as follows.  In Section \ref{sec:model}, we introduce our generic solid mechanics problem formulation, as well as a simplified model problem based on these equations: the linear elastic wave propagation problem.  Sections \ref{sec:schwarz_fom} and \ref{sec:schwarz_opinf} focus on NO-SAM for FOM coupling and OpInf coupling, respectively.  We first provide the formulation of the method with two types of transmission BCs: alternating Dirichlet-Neumann and Robin-Robin. The resulting DN-NO-SAM and RR-NO-SAM methods are then evaluated on our model problem in terms of both accuracy and efficiency. For all the couplings considered (FOM-FOM, FOM-OpInf, OpInf-OpInf), RR-NO-SAM with the optimal selection of free parameters typically delivers the most accurate and efficient solution. We discuss conclusions and future work in Section \ref{sec:conc}.

\section{Solid mechanics problem formulation} \label{sec:model}

Consider the strong form Euler-Lagrange equations for a general dynamic solid mechanics problem:
\begin{equation}
    \label{eqn:dynamic_elasticity_pde}
    \mathrm{Div}\, ~P + \rho_0 ~B = \rho_0 \ddot{\boldsymbol{\Psi}} \quad \text{in } \Omega \times I
\end{equation}

\noindent where $\Omega \in \mathbb{R}^d$ for $d \in \{1, 2, 3\}$ is an open bounded domain, $I \coloneqq \{t \in [t_0, t_f]\}$ is a closed time interval with $0 \leq t_0<t_f$, and $~x=\boldsymbol{\Psi}(~X, t):\Omega \times I \rightarrow \mathbb{R}^d$ describes the mapping from the reference configuration to the current configuration with $~X \in \Omega$. The overdot notation in $\ddot{\boldsymbol{\Psi}}$ denotes a time derivative such that $\dot{\boldsymbol{\Psi}} \coloneqq \frac{\partial \boldsymbol{\Psi}}{\partial t}$ and $\ddot{\boldsymbol{\Psi}} \coloneqq \frac{\partial^2 \boldsymbol{\Psi}}{\partial t^2}$. The symbol $~P$ denotes the first Piola-Kirchoff stress, $~B$ is the specific body force, and $\rho_0$ is the mass density defined in the reference configuration. It is important to note that the constitutive model used to describe $~P$ is general and can include nonlinear models.

We then have the following initial and boundary conditions for the governing equation described in \eqref{eqn:dynamic_elasticity_pde}:
\begin{equation}
    \label{eqn:dynamic_elasticity_bcs}
    \begin{aligned}
    \boldsymbol{\Psi}(~X, t_0) &= ~X_0 && \text{in } \Omega \\
    \dot{\boldsymbol{\Psi}}(~X, t_0) &= ~v_0 && \text{in } \Omega \\
    \boldsymbol{\Psi}(~X, t) &= \boldsymbol{\chi} && \text{on } \partial_{\Psi} \Omega \times I \\
    ~P ~N &= ~T && \text{on } \partial_{T} \Omega \times I
    \end{aligned}
\end{equation}

\noindent where $\partial_{\Psi} \Omega$ and $\partial_{T} \Omega$ denote the Dirichlet and Neumann portions, respectively, of the domain boundary $\partial \Omega$ with $\partial \Omega = \partial_{\Psi} \Omega \cup \partial_T \Omega$ and $\partial_{\Psi} \Omega \cap \partial_T \Omega = \emptyset$, and $~N$ denotes the unit vector normal to $\partial \Omega$. 

The weak variational form of \eqref{eqn:dynamic_elasticity_pde} is then given by:


\begin{equation}
    \label{eqn:weak_form}
    \int_I \left[ \int_{\Omega} \left( \rho_0 \ddot{\boldsymbol{\Psi}}\cdot \boldsymbol{\xi} + ~P : \mathrm{Grad} \, ~\xi - \rho_0 ~B\cdot \boldsymbol{\xi} \right) \, d \Omega  - \int_{\partial_T \Omega} \left( ~P ~N \right)\cdot \boldsymbol{\xi}  \, d S  \right] \, dt = 0
\end{equation}

\noindent where $\boldsymbol{\xi}$ is a test function defined so that $\mathcal{V} \coloneqq \{ \boldsymbol{\xi} \in W_2^1(\Omega \times I): \boldsymbol{\xi=~0} \text{ on } \partial_{\Psi} \Omega \times I \ \cup \Omega \times t_0 \cup \Omega \times t_f\}$. 
The variational form in \eqref{eqn:weak_form} is discretized in space using the classical Galerkin finite element method (FEM) \cite{hughes1987finite}.  The Neumann boundary condition from \eqref{eqn:dynamic_elasticity_bcs} is implemented by inserting $~P~N=~T$ into the second integral of \eqref{eqn:weak_form}. 

For the sake of concreteness, we consider in this work a specific linear elastic initial boundary value problem (IBVP) for the equations \eqref{eqn:dynamic_elasticity_pde}--\eqref{eqn:dynamic_elasticity_bcs} known as the one-dimensional (1D) linear elastic wave propagation problem, described in detail in Section \ref{sec:problem_setting}.
It is straightforward to show that, in the case of a linear elastic material, the semi-discrete matrix-vector form of \eqref{eqn:weak_form} can be written as:
\begin{equation}
    \label{eqn:matrix_vector_form}
    ~M \ddot{~u}+ ~K ~u = ~F
\end{equation}

\noindent where $~M \in \mathbb{R}^{N \times N}$ is the (consistent) mass matrix, $~K \in \mathbb{R}^{N \times N}$ is the stiffness matrix, and $~F \in \mathbb{R}^{N}$ is the vector of applied external forces. Here, $~u \in \mathbb{R}^{N}$ is the vector of displacement degrees of freedom, defined as $~u := \boldsymbol{\Psi}(~X,t) - ~X$, and $\ddot{~u}$ denotes the acceleration.

\subsection{Model problem: wave propagation in a 1D linear elastic bar} \label{sec:problem_setting}


In order to explore the various transmission conditions available when formulating NO-SAM, we consider a 1D linear elastic wave propagation problem based on the governing 
equations \eqref{eqn:dynamic_elasticity_pde}. This problem, considered previously in \cite{Mota:2022} and \cite{Barnett:2022Schwarz}, is chosen due to its ability to reveal coupling artifacts, such as spurious oscillations, as the wave passes through the interface between subdomains. 

Consider a simple bar geometry of length $L$ with $\Omega = [0,L] \in \mathbb{R}$. Both ends of the bar are clamped such that $u(0, t)=u(L, t)=0$ for all $t \geq t_0$. The initial conditions are that of a symmetric Gaussian wave described by the following equation:
\begin{equation}
    \label{eqn:gaussian_wave}
    u(x, t_0)=a \, \, \mathrm{exp} \left[-\frac{(x-b)^2}{2w^2} \right]
\end{equation}

\noindent for $a, b, w \in \mathbb{R}$, where these parameters represent the peak height, center, and standard deviation of the Gaussian wave, respectively.

To benchmark the performance of each transmission condition for NO-SAM, we developed a Python script that implements the numerical framework described above. All numerical examples in this work use the following problem setting, unless otherwise specified. A linear elastic constitutive model is employed with material constants $E=1 \, GPa$ and $\rho =1000 \, kg \, m^{-3}$. The bar is of length $L=1 \, m$ and is discretized into uniform elements of length $h=0.001 \, m$. The domain decomposition of the bar is defined as $\Omega_1 = [0, 0.6]$ and $\Omega_2 = [0.6, 1]$. The semi-discrete equation \eqref{eqn:matrix_vector_form} is advanced in time using the Newmark-$\beta$ time integration scheme of constant average acceleration, for which $\beta=0.25$ and $\gamma=0.5$. The simulation is run from initial time $t_0=0 \, s$ to final time $t_f=10^{-3} \, s$ using a time step of $\Delta t=2.5 \times 10^{-7} \, s$. This time step was chosen based on the wave propagation speed $c=\sqrt{\frac{E}{\rho}}=1000 \, m \, s^{-1}$, for which the stable time step was calculated to be $\Delta t_{stable}=\frac{h}{c}=10^{-6} \, s$ \cite{Mota:2022, Barnett:2022Schwarz}. Uniform, conformal meshes and identical time-step sizes are used for each subdomain, although SAM does not require this restriction. The initial displacement is defined using the symmetric Gaussian wave in \eqref{eqn:gaussian_wave} with $a=0.005\, m$, $b=\frac{L}{2}=0.5\, m$, and $w=0.02\, m$. 

\section{NO-SAM for FOM coupling} \label{sec:schwarz_fom}

This section outlines the Schwarz alternating method for the domain-decomposition based coupling between multiple non-overlapping subdomains. The formulation is presented for a general setting to accommodate different types of transmission conditions between the subdomains.

\subsection{NO-SAM formulation} \label{sec:schwarz_fom_formulation}

\begin{figure}[htbp!]
    \centering
    \includegraphics[width=0.5\linewidth]{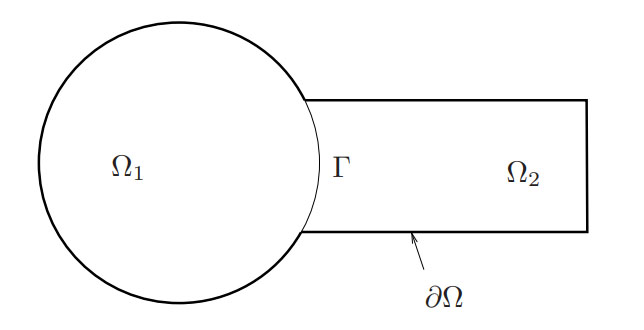}
    \caption{Illustration showing an example non-overlapping domain decomposition.}
    \label{fig:non_overlapping_decomp}
\end{figure}

Without loss of generality, consider a partition of a given domain $\Omega$ into two non-overlapping subdomains, $\Omega_1$ and $\Omega_2$, as shown in Figure \ref{fig:non_overlapping_decomp}. Here, \textit{non-overlapping} indicates that the subdomains satisfy $\Omega_1 \cap \Omega_2 = \emptyset$. We are interested in solving the weak variational equation of the form \eqref{eqn:weak_form} using this decomposition. 

As described in more detail in \cite{Mota:2022}, the dynamic SAM algorithm operates by discretizing the temporal domain into time intervals $I_n := \left[ t_n, t_{n+1} \right]$ and solving a sequence of subdomain-local PDEs to convergence within each time interval.  Specifically, for a given time interval $I_n$, NO-SAM solves the following sequence of subdomain-local problems:
\begin{equation}
    \label{eqn:schwarz_1}
    \left\{\begin{array}{rcll}
    ~M_1 \ddot{~u}_1^{(s+1)} + ~K_1 ~u_1^{(s+1)} &=&~F_1^{(s+1)}& \text{in } \Omega_1
    \\[1.0ex]
    ~u_1^{(s+1)} &=& ~u_D 
    & \text{on } \partial_{\Psi} \Omega_1 \backslash \Gamma
    \\[0.5ex]
    \alpha_{12} ~T_1^{(s+1)} + \beta_{12} ~u_1^{(s+1)} &=& ~\lambda^{(s+1)}_1 
       & \text{on } \Gamma
    \end{array}\right.
\end{equation}

\begin{equation}
    \label{eqn:schwarz_2}
    \left\{\begin{array}{rcll}
    ~M_2 \ddot{~u}_2^{(s+1)} + ~K_2 ~u_2^{(s+1)} &=&~F_2^{(s+1)}& \text{in } \Omega_2
    \\[1.0ex]
    ~u_2^{(s+1)} &=& ~u_D 
    & \text{on } \partial_{\Psi} \Omega_2 \backslash \Gamma
    \\[0.5ex]
    \alpha_{21} ~T_2^{(s+1)} + \beta_{21} ~u_2^{(s+1)} &=& ~\lambda^{(s+1)}_2 
       & \text{on } \Gamma
    \end{array}\right.
\end{equation}

\begin{align}
\label{eqn:relaxation_1}
~\lambda^{(s+1)}_1 &= \theta_1 \left[\alpha_{12} ~T_2^{(s)} + \beta_{12} ~u_2^{(s)}\right]
+ (1-\theta_1)\,~\lambda^{(s)}_1
\\[8pt]
\label{eqn:relaxation_2}
~\lambda^{(s+1)}_2 &= \theta_2 \left[\alpha_{21} ~T_1^{(s+1)} + \beta_{21} ~u_1^{(s+1)}\right]
+ (1-\theta_2)\,~\lambda^{(s)}_2
\end{align}

\noindent where $s$ denotes the Schwarz iteration index and $\theta_1, \theta_2 \in (0, 1]$ are relaxation parameters. For convenience, we rewrite the Dirichlet boundary condition defined in \eqref{eqn:dynamic_elasticity_bcs} in terms of displacement, denoting the boundary condition on the outer domain boundary as $~u_D$. At each iteration, the displacements $~u^{(s+1)}$ are solved for, and the solution on the interface $\Gamma$ is used to define the transmission conditions for the neighboring subdomain. To allow for a general choice of transmission conditions, we introduce constants $\alpha_{ij}, \beta_{ij} \in \mathbb{R}$ for $i,j = 1, \dots, N_k$, where $N_k$ is the total number of subdomains ($N_k=2$ in the present study). The reader can observe that, by appropriate selection of these constants, Dirichlet, Neumann, and Robin transmission conditions can all be represented in \eqref{eqn:schwarz_1}--\eqref{eqn:relaxation_2}. 

The iteration defined in \eqref{eqn:schwarz_1}--\eqref{eqn:relaxation_2} is repeated until convergence. We adopt the following relative convergence criterion 
\begin{equation}
\label{eqn:schwarz_tolerance}
  \frac{\big\lVert \left[ ~u^{(s+1)} - ~u^{(s)} \right] 
  + \Delta t\, \left[ ~v^{(s+1)} - ~v^{(s)} \right]
  + \tfrac{1}{2}\,\Delta t^2\, \left[ ~a^{(s+1)} - ~a^{(s)} \right]\big\rVert_2}{
  \big\lVert  ~u^{(s)} 
  + \Delta t\, ~v^{(s)}
  + \tfrac{1}{2}\,\Delta t^2\, ~a^{(s)} \big\rVert_2} 
  < \delta
\end{equation}

\noindent where $\delta$ denotes the prescribed convergence tolerance. This criterion requires that the tolerance be satisfied across all subdomains before declaring convergence of the Schwarz iteration, thereby ensuring that all displacement, velocity, and acceleration solutions have converged prior to advancing to the next time step. In all numerical examples presented here, we use $\delta = 10^{-8}$.

This work focuses specifically on the alternating Dirichlet-Neumann and Robin-Robin transmission conditions, for which the $\alpha_{ij}$ and $\beta_{ij}$ constants are defined as

\vspace{5pt} 
\begin{itemize}
    \setlength{\itemsep}{5pt} 
    \item \textbf{Alternating Dirichlet-Neumann:} $\alpha_{12} = \beta_{21} = 0$, and $\alpha_{21}, \beta_{12} = 1$
    \item \textbf{Robin--Robin:} $\alpha_{12}, \alpha_{21}, \beta_{12}, \beta_{21} \neq 0$
\end{itemize}
\vspace{5pt} 
We consider formulations with these specific transmission conditions, as they have shown to lead to a convergent NO-SAM iteration in \cite{Zanolli:1987} and \cite{Lions:1990}, respectively. Relaxation is typically introduced only to subdomains using a Dirichlet transmission condition on the boundary $\Gamma$, as occurs in the alternating Dirichlet-Neummann variant of NO-SAM, where $\alpha_{ij} = 0$ and $\beta_{ij} \neq 0$. 
It has been demonstrated in the literature that, for certain PDEs, setting $\theta_i<1$ can reduce the number of Schwarz iterations required for convergence \cite{Zanolli:1987}.  In this work, $~\lambda^i_0$ is initialized to zero for the first iteration, i.e.,  $~\lambda^i_0 = ~0$.

For the case of the Robin-Robin NO-SAM variant, the Robin transmission condition is applied as follows. Consider without loss of generality the subdomain problem in $\Omega_1$.  In this case, the third line of \eqref{eqn:schwarz_1} is rearranged to solve for $~T_1^{(s+1)}$, and the resulting expression is then substituted into the external force term $~F_1^{(s+1)}$ of the first line of \eqref{eqn:schwarz_1}, yielding
\begin{equation}
    \label{eqn:robin_final}
    ~M_1 \ddot{~u}_1^{(s+1)} + \left[~K_1 + \frac{\beta_{12}}{\alpha_{12}}~S_1\right] ~u_1^{(s+1)} =~F_1^{(s+1)} + \frac{1}{\alpha_{12}}~R_1^{(s+1)}
\end{equation}



\noindent where the additional $~S$ and $~R$ terms can be viewed as Robin contributions to the stiffness matrix and external force vector, respectively. Their derivation is provided in the Appendix. The corresponding derivation for $\Omega_2$ is analogous and omitted for brevity.
 \\

\noindent {\bf Remark 2.1. }There exists a wide body of literature focused on how the parameters $\alpha_{ij}$, $\beta_{ij}$ and $\theta_i$ in \eqref{eqn:schwarz_1}--\eqref{eqn:schwarz_2} should be selected to get the best convergence of SAM for a given IBVP.  The resulting methods are termed ``Optimized Schwarz methods"; for a succinct overview, the reader is referred to \cite{Gander:2006}.  While it is possible to derive optimal values of the relevant parameters analytically for certain simple PDEs, this is not generally the case for generic nonlinear PDEs based on \eqref{eqn:weak_form}.  Herein, we find the most viable values of  $\alpha_{ij}$, $\beta_{ij}$ and $\theta_i$ through a numerical parameter sweep study.\\

\subsection{Numerical results} \label{sec:schwarz_fom_results}

Our first set of numerical results examines the influence of the Robin boundary condition parameters, $\alpha_{ij}$ and $\beta_{ij}$, as defined in Section \ref{sec:schwarz_fom_formulation}, on: (i) the error between the converged Schwarz solution and the corresponding monolithic solution, and (ii) the number of Schwarz iterations per time step required to satisfy the prescribed convergence tolerance. This is achieved by performing a parameter sweep for all four $\alpha_{ij}$ and $\beta_{ij}$ parameters within the 1D wave propagation problem setting described in Section \ref{sec:problem_setting}. The parameter spaces used in the sweep are:
\begin{equation}
\label{eqn:parameter_space}
\begin{aligned}
    \alpha_{ij} &\in  
        \{ 10^{-3},\, 10^{-1},\, 1,\, 3,\, 5 \} \cdot\tfrac{1}{\sigma_{\max}} \\
    \beta_{ij}  &\in 
        \{ 10^{-3},\, 10^{-1},\, 1,\, 3,\, 5 \}
\end{aligned}
\end{equation}

\noindent where $\sigma_{\max} := \max_{t,\,x} \sigma(x,t)$ denotes the maximum stress within the bar for the corresponding monolithic solution. Normalizing the $\alpha_{ij}$ parameter by a reference stress value compensates for the large difference in orders of magnitude between the displacements and tractions used in the Robin transmission conditions and improves the interpretability of the results. The parameter spaces defined in \eqref{eqn:parameter_space} were selected to capture the full range of variation in both the error and number of iterations required for convergence. 

The error between the converged Schwarz solution and the corresponding monolithic solution is calculated at each time step. For subdomain $k$ at time $t_n$ we define
\begin{equation}
  \label{eqn:error_individual}
  \varepsilon_k(t_n) := \big\lVert \Delta~u_k(t_n)\big\rVert_2
  + \Delta t\,\big\lVert \Delta~v_k(t_n)\big\rVert_2
  + \tfrac{1}{2}\,\Delta t^2\,\big\lVert \Delta~a_k(t_n)\big\rVert_2
\end{equation}

\noindent where
\begin{equation}
  \label{eqn:vector_differences}
  \begin{aligned}
    \Delta~u_k(t_n) &= ~u_k^{\mathrm{SCH}}(t_n) - ~u_k^{\mathrm{MON}}(t_n) \\
    \Delta~v_k(t_n) &= ~v_k^{\mathrm{SCH}}(t_n) - ~v_k^{\mathrm{MON}}(t_n) \\
    \Delta~a_k(t_n) &= ~a_k^{\mathrm{SCH}}(t_n) - ~a_k^{\mathrm{MON}}(t_n)
  \end{aligned}
\end{equation}

\noindent and \(\lVert\cdot\rVert_2\) denotes the Euclidean norm. Here \(~u^{\mathrm{SCH}}\), \(~v^{\mathrm{SCH}}\), \(~a^{\mathrm{SCH}}\) are the converged Schwarz displacement, velocity and acceleration, and \(~u^{\mathrm{MON}}\), \(~v^{\mathrm{MON}}\), \(~a^{\mathrm{MON}}\) are the corresponding monolithic quantities. Lastly, the time-averaged error is defined as
\begin{equation}
    \label{eqn:error_avg}
        \varepsilon_{avg} :=\frac{1}{N_t-1} \sum_{k=1}^{N_k} \sum_{n=2}^{N_t} \varepsilon_k(t_n)
\end{equation}

\noindent where $N_k$ denotes the number of subdomains, and $N_t$ is the number of time-steps. 

\begin{figure}[b!]
    \begin{center}
        {\includegraphics[width=0.95\textwidth]
        {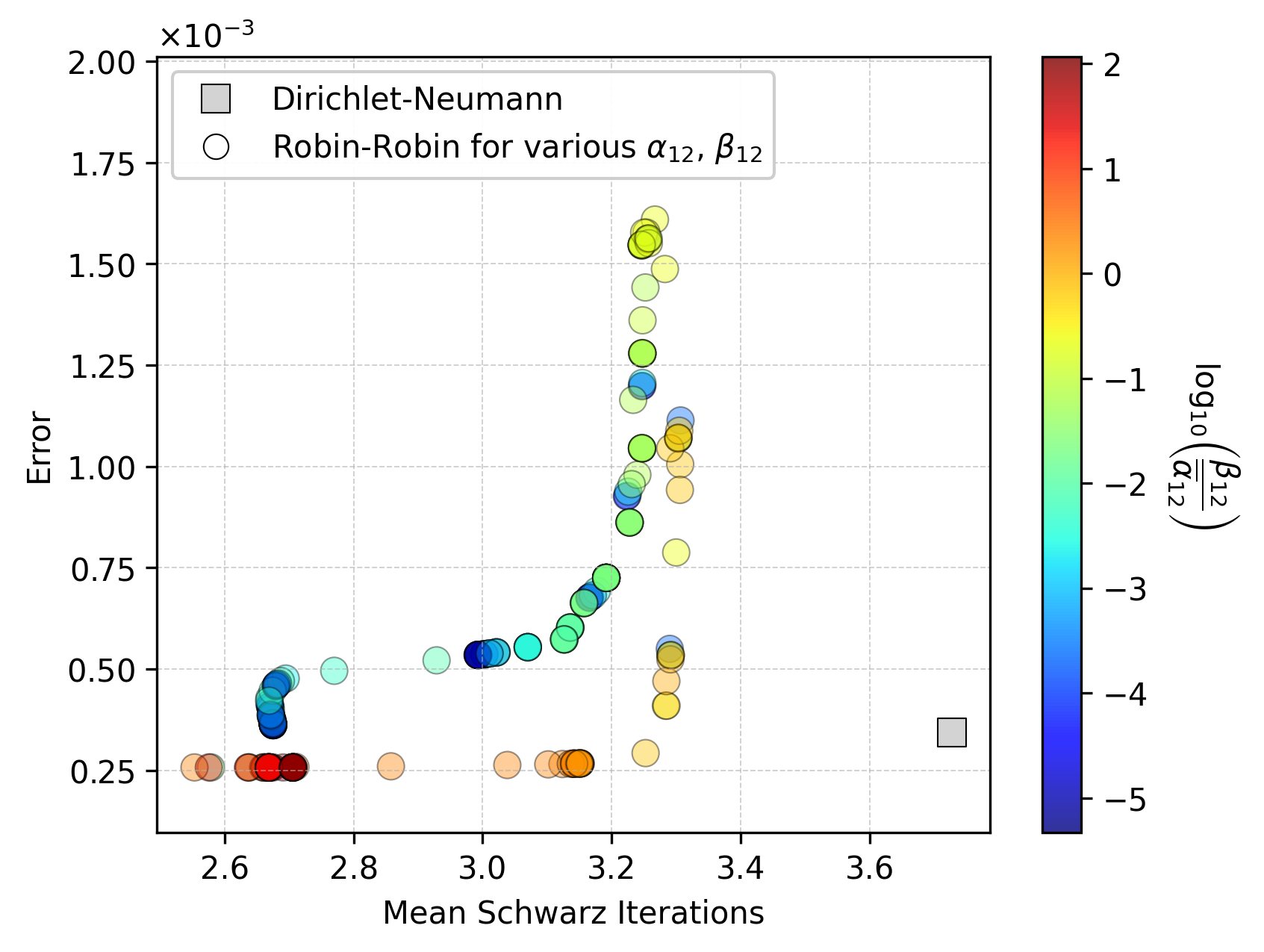}}
    \end{center}
    \caption{FOM-FOM coupling: Pareto plot illustrating the relationship between the error in the converged Schwarz solution relative to the corresponding monolithic solution and the mean number of Schwarz iterations per time step required for convergence, for varying $\alpha_{ij}$ and $\beta_{ij}$ parameters. The unrelaxed alternating Dirichlet-Neumann result is also plotted as a reference in gray.}
    \label{fig:FOM_FOM_robin_pareto}
\end{figure}

A Pareto plot of the parameter sweep results is presented in Figure \ref{fig:FOM_FOM_robin_pareto}. The error, as defined in~\eqref{eqn:error_avg}, is plotted against the mean number of Schwarz iterations per time step required for convergence for varying $\alpha_{ij}$ and $\beta_{ij}$ parameters. The color mapping represents the ratio $\beta_{12} / \overline{\alpha}_{12}$, serving as an indicator of the relative magnitude of Dirichlet to Neumann data used in the Robin transmission condition. Here, $\overline{\alpha}_{ij} := \sigma_{max}\,\cdot \alpha_{ij}$. A high ratio (red) approaches a pure Dirichlet condition, whereas a low ratio (blue) approaches a pure Neumann condition.  
It is curious that the near Dirichlet-Dirichlet case gives the most accurate and efficient results, given that Dirichlet-Dirichlet NO-SAM is known not to converge \cite{Lions:1990}. This behavior was also observed in our implementation.
For brevity, the plot showing results for $\beta_{21} / \overline{\alpha}_{21}$ is omitted, as similar trends are observed. The result obtained with an alternating Dirichlet–Neumann transmission condition is also included in gray as a reference.  No relaxation was employed in generating our alternating Dirichlet-Neumann results, as numerical experiments revealed that the inclusion of relaxation worsened both accuracy and efficiency of the NO-SAM algorithm.

The results summarized herein indicate that parameter choices near the relative extremes of $\alpha_{ij}$ and $\beta_{ij}$ yield both a lower error and fewer Schwarz iterations to reach convergence. Moving away from these extremes toward intermediate values, where the Dirichlet and Neumann data are of comparable magnitude, results in a significant increase in both quantities. The outer contour of the plot illustrates this transition across the parameter range. The key cases are summarized in Table~\ref{tab:pareto_results}. Notably, all Robin parameter sets within the explored space converged in fewer iterations than the reference Dirichlet–Neumann result, although errors vary significantly. Additional $\alpha_{ij}$ and $\beta_{ij}$ parameters outside the space defined in \eqref{eqn:parameter_space} were generally found to fall within the contour outlined in Figure~\ref{fig:FOM_FOM_robin_pareto}, but were omitted for the sake of clarity.

\begin{table}[th!]
    \centering
    \caption{FOM-FOM coupling: summary of parameter sets and their corresponding error and iteration counts.}
    \label{tab:pareto_results}
    \begin{tabular}{lcccccc}
        \toprule[1pt]
        Case & $\overline{\alpha}_{12}$ & $\overline{\alpha}_{21}$ & $\beta_{12}$ & $\beta_{21}$ & Error & Iterations\\
        \midrule[1pt]
        Lowest Error  & $10^{-3}$ & $10^{-3}$ & $1$ & $1$ & $2.57 \times 10^{-4}$ & $2.66$\\
        Highest Error & $10^{-1}$ & $1$ & $1$ & $3$ & $1.61 \times 10^{-3}$ & $3.27$\\
        Lowest Iterations & $10^{-3}$ & $10^{-3}$ & $10^{-1}$ & $5$ & $2.57 \times 10^{-4}$ & $2.55$\\
        Highest Iterations & $10^{-1}$ & $10^{-1}$ & $10^{-3}$ & $3$ & $1.11 \times 10^{-3}$ & $3.31$\\
        Dirichlet-Neumann & -- & -- & -- & -- & $3.43\times 10^{-4}$ & $3.73$ \\
        \bottomrule[1pt]
    \end{tabular}
\end{table}

\begin{figure}[h!]
    \centering
    
    \begin{subfigure}{0.45\linewidth}
        \centering
        \includegraphics[width=\linewidth]{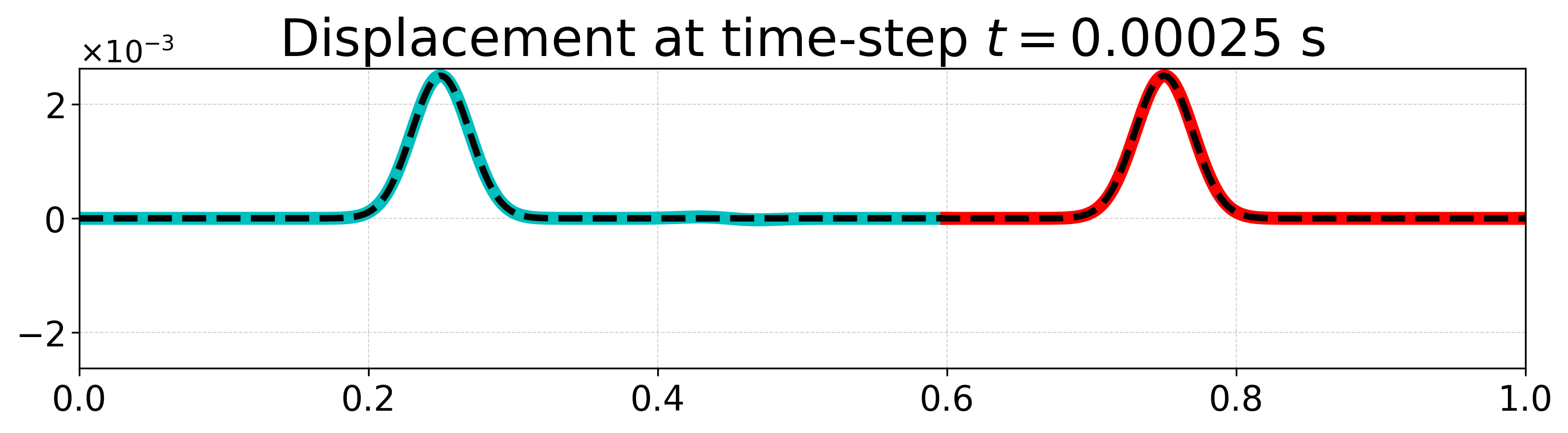}
    \end{subfigure}
    \hspace{3mm}
    \begin{subfigure}{0.45\linewidth}
        \centering
        \includegraphics[width=\linewidth]{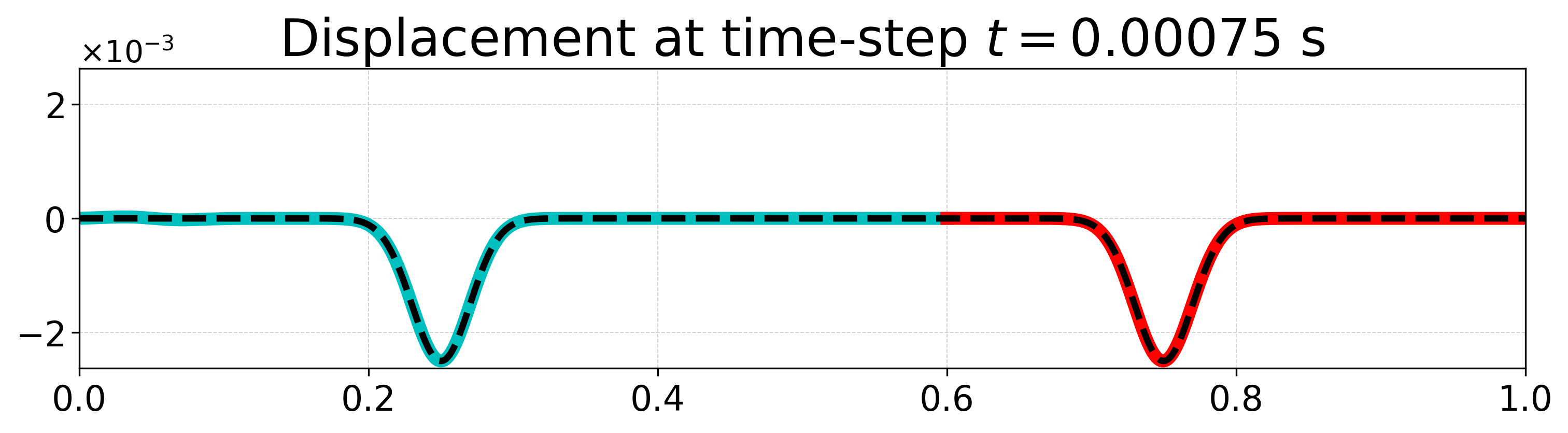}
    \end{subfigure}
    
    \vspace{1mm}
    
    \begin{subfigure}{0.45\linewidth}
        \centering
        \includegraphics[width=\linewidth]{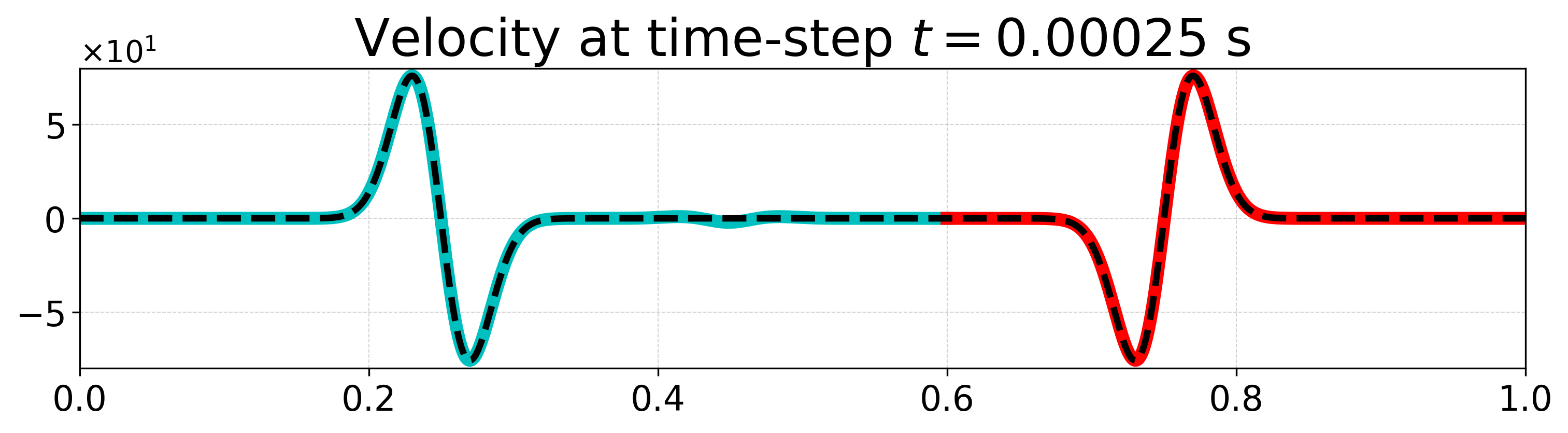}
    \end{subfigure}
    \hspace{3mm}
    \begin{subfigure}{0.45\linewidth}
        \centering
        \includegraphics[width=\linewidth]{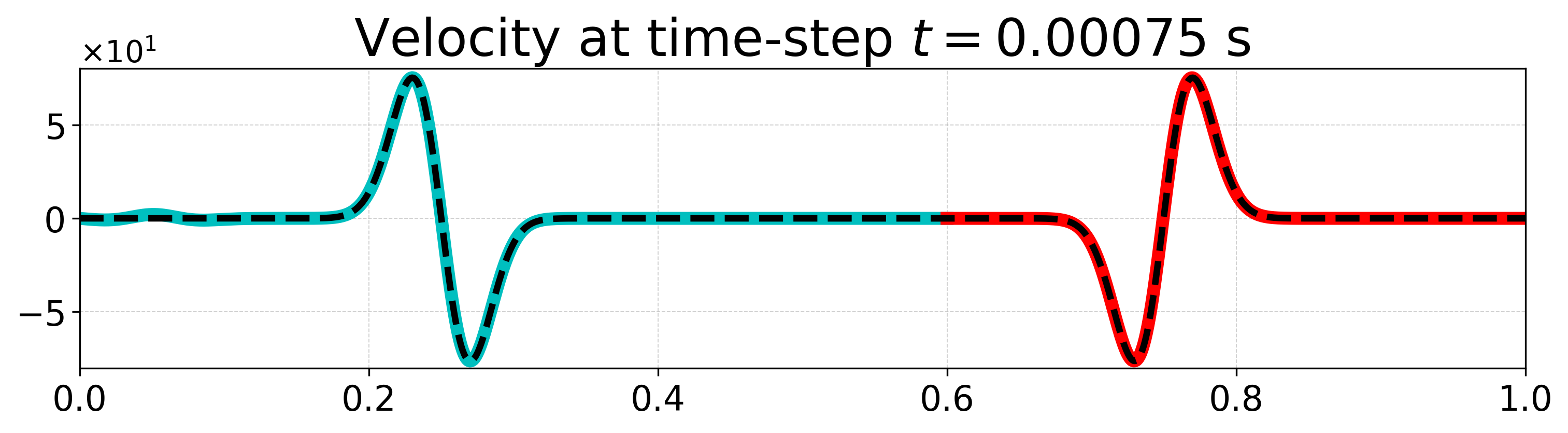}
    \end{subfigure}
    
    \vspace{1mm}
    
    \begin{subfigure}{0.45\linewidth}
        \centering
        \includegraphics[width=\linewidth]{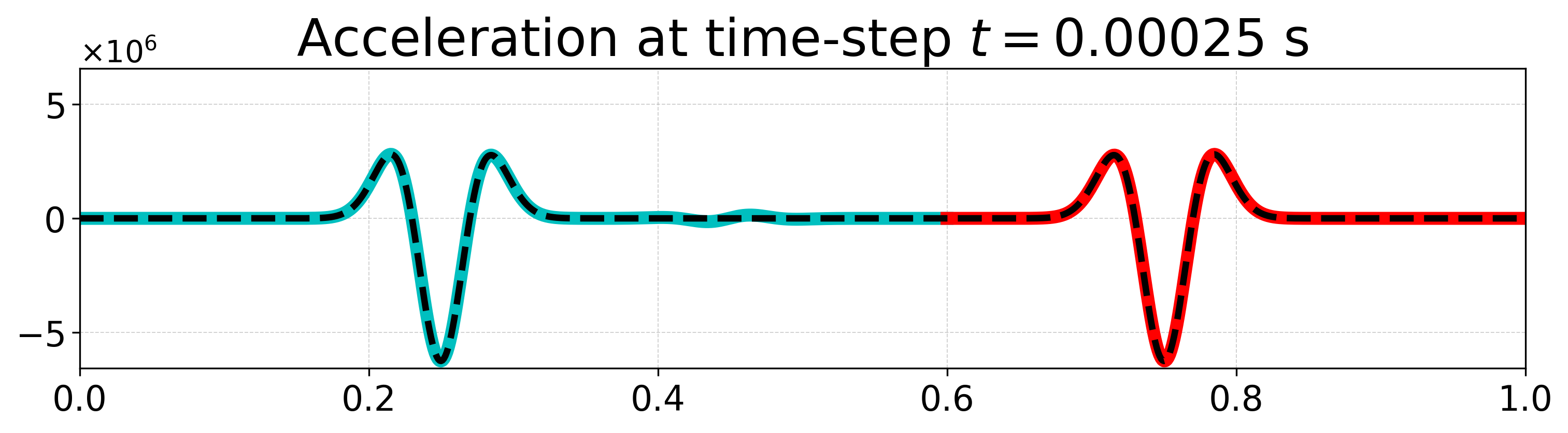}
        \caption{$t=2.5 \times 10^{-4} \, s$}
    \end{subfigure}
    \hspace{3mm}
    \begin{subfigure}{0.45\linewidth}
        \centering
        \includegraphics[width=\linewidth]{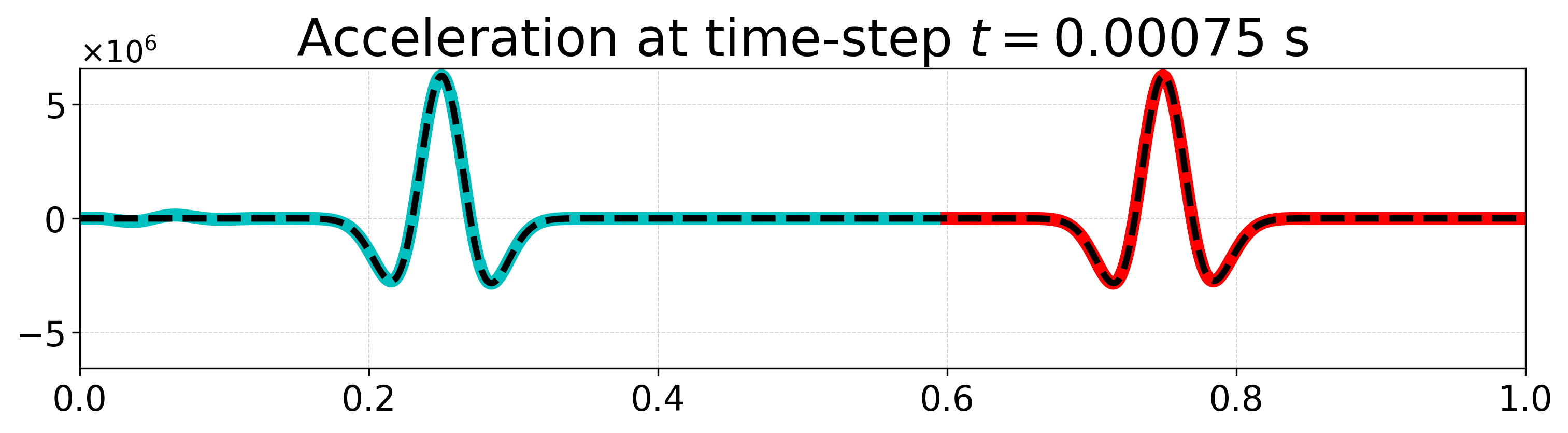}
        \caption{$t=7.5 \times 10^{-4} \, s$}
    \end{subfigure}
    
    \caption{FOM-FOM coupling: plots of the NO-SAM displacement, velocity, and acceleration solutions obtained with alternating Dirichlet–Neumann transmission conditions, compared to the corresponding monolithic solution (dashed line). The blue subdomain $\Omega_1$ employs a Dirichlet transmission condition, and the red subdomain $\Omega_2$ employs a Neumann transmission condition.}
    \label{fig:DN_plots}
\end{figure}

\begin{figure}[h!]
    \centering
    
    \begin{subfigure}{0.45\linewidth}
        \centering
        \includegraphics[width=\linewidth]{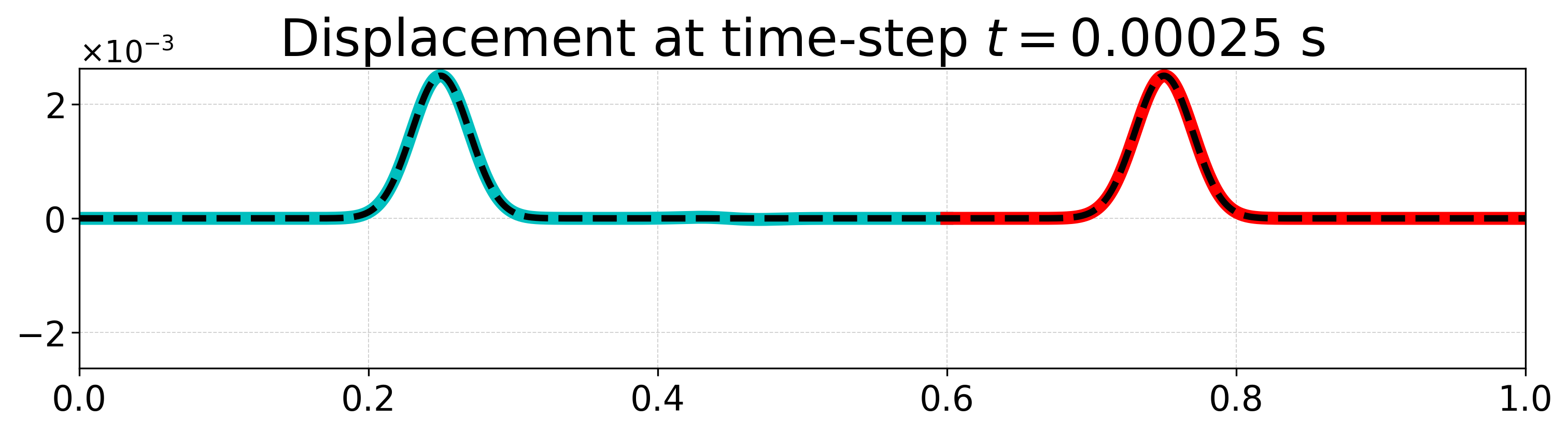}
    \end{subfigure}
    \hspace{3mm}
    \begin{subfigure}{0.45\linewidth}
        \centering
        \includegraphics[width=\linewidth]{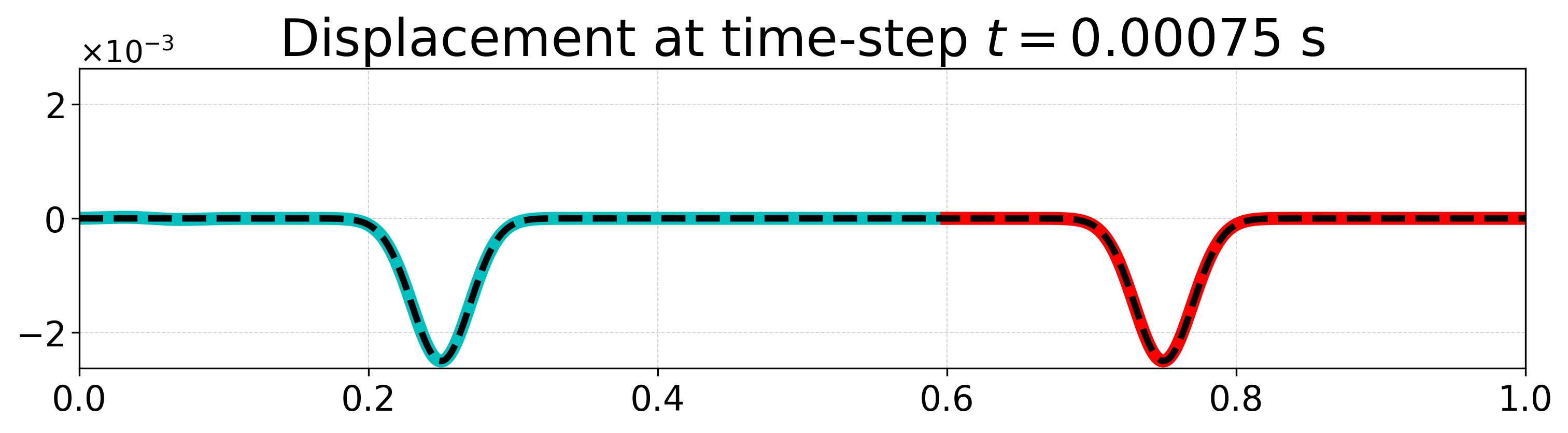}
    \end{subfigure}
    
    \vspace{1mm}
    
    \begin{subfigure}{0.45\linewidth}
        \centering
        \includegraphics[width=\linewidth]{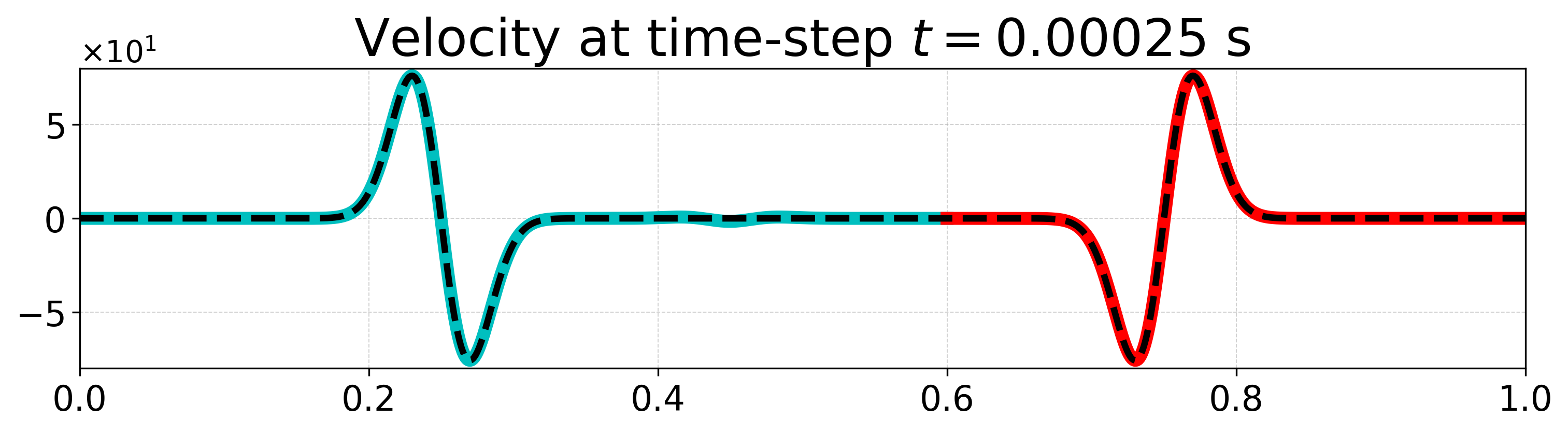}
    \end{subfigure}
    \hspace{3mm}
    \begin{subfigure}{0.45\linewidth}
        \centering
        \includegraphics[width=\linewidth]{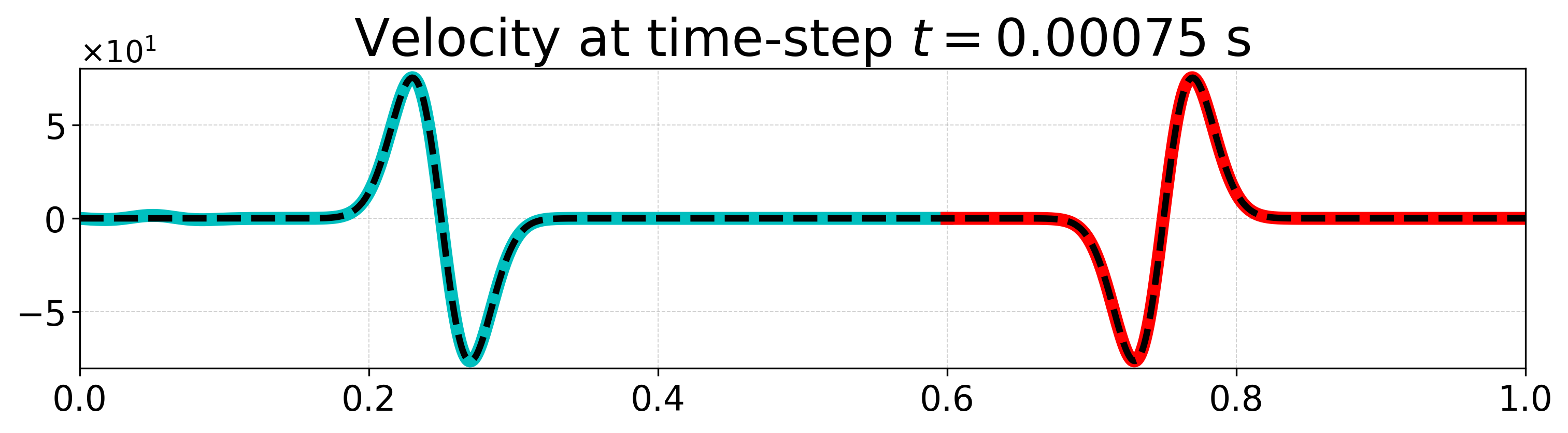}
    \end{subfigure}
    
    \vspace{1mm}
    
    \begin{subfigure}{0.45\linewidth}
        \centering
        \includegraphics[width=\linewidth]{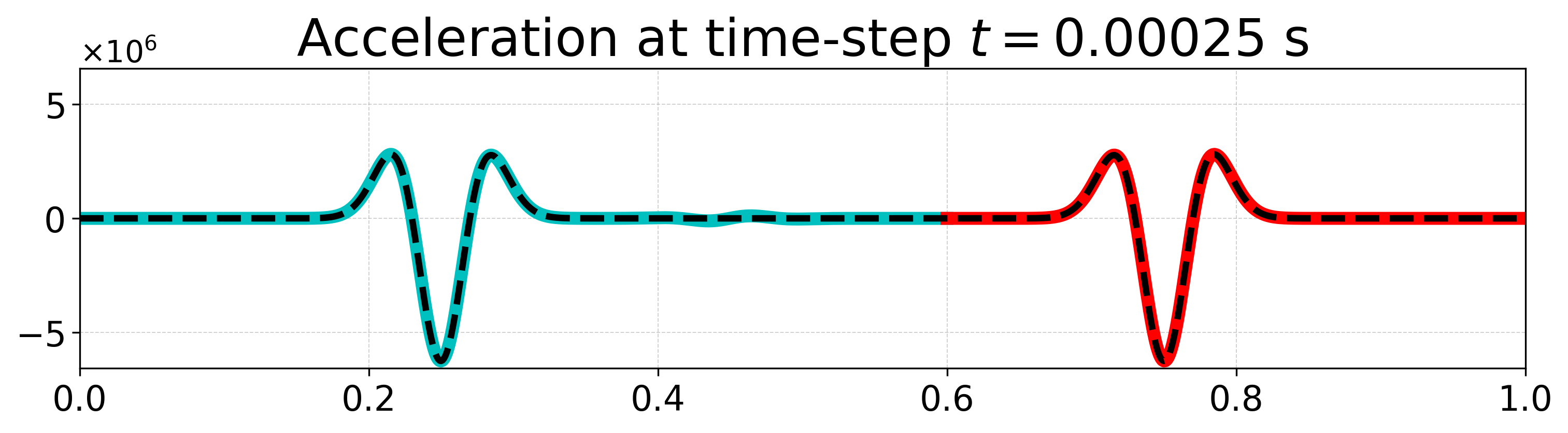}
        \caption{$t=2.5 \times 10^{-4} \, s$}
    \end{subfigure}
    \hspace{3mm}
    \begin{subfigure}{0.45\linewidth}
        \centering
        \includegraphics[width=\linewidth]{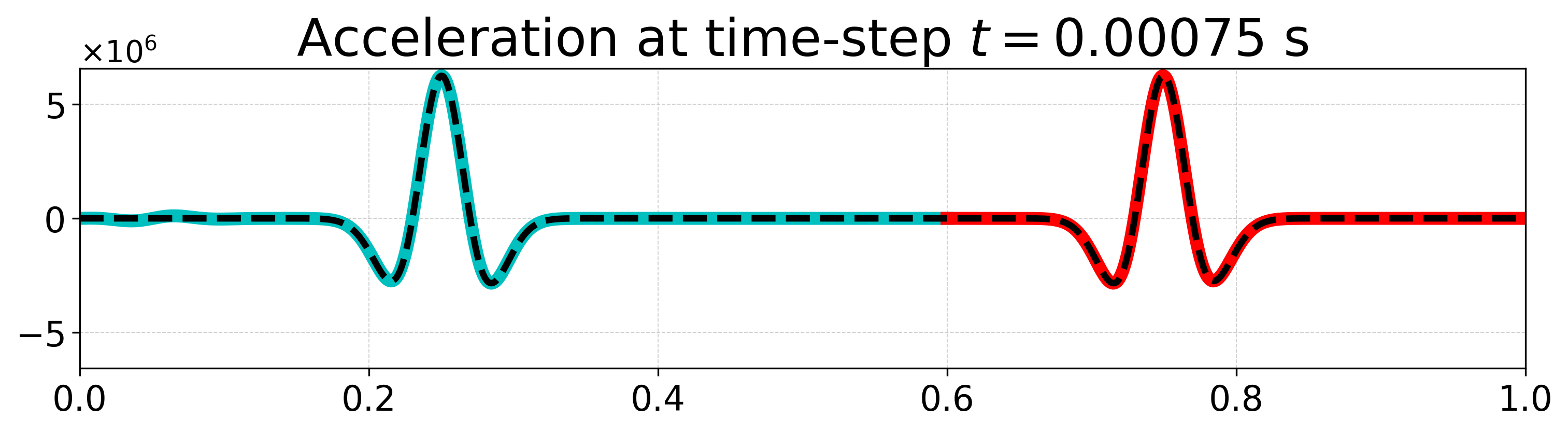}
        \caption{$t=7.5 \times 10^{-4} \, s$}
    \end{subfigure}
    \caption{FOM-FOM coupling: plots of the NO-SAM displacement, velocity, and acceleration solutions obtained with Robin–Robin transmission conditions (lowest error parameter set, $\left[ \overline{\alpha}_{12}, \overline{\alpha}_{21}, \beta_{12},  \beta_{21} \right] = \left[ 10^{-3}, 10^{-3}, 1, 1\right]$) compared to the corresponding monolithic solution (dashed line).  The solutions computed for subdomains $\Omega_1$ and $\Omega_2$ are shown in blue and red, respectively.
}
    \label{fig:RR_Best_plots}
\end{figure}

\begin{figure}[h!]
    \centering
    
    \begin{subfigure}{0.45\linewidth}
        \centering
        \includegraphics[width=\linewidth]{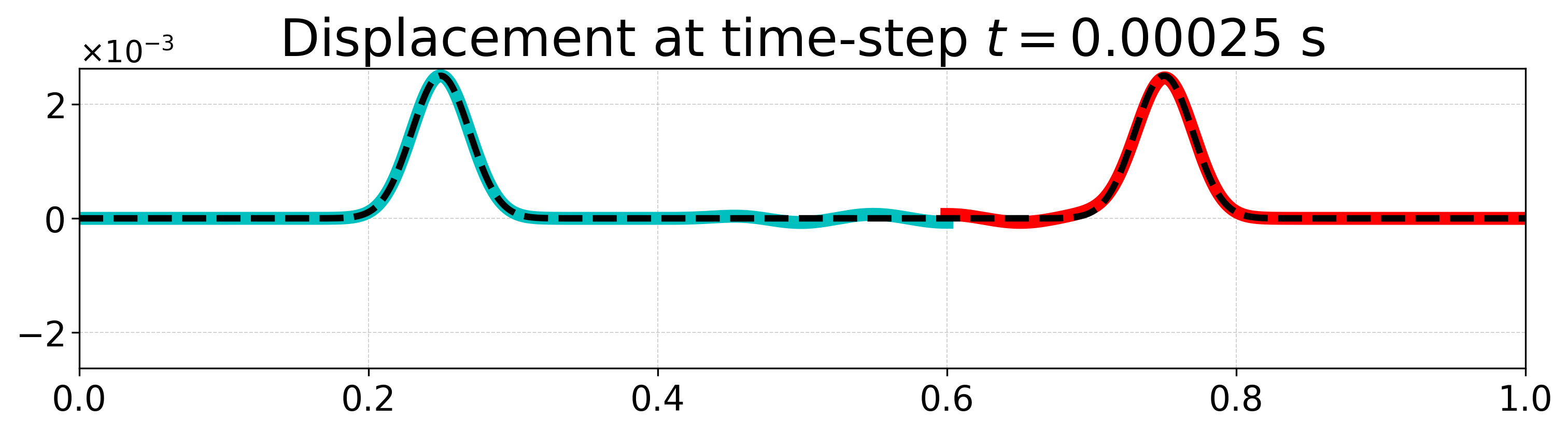}
    \end{subfigure}
    \hspace{3mm}
    \begin{subfigure}{0.45\linewidth}
        \centering
        \includegraphics[width=\linewidth]{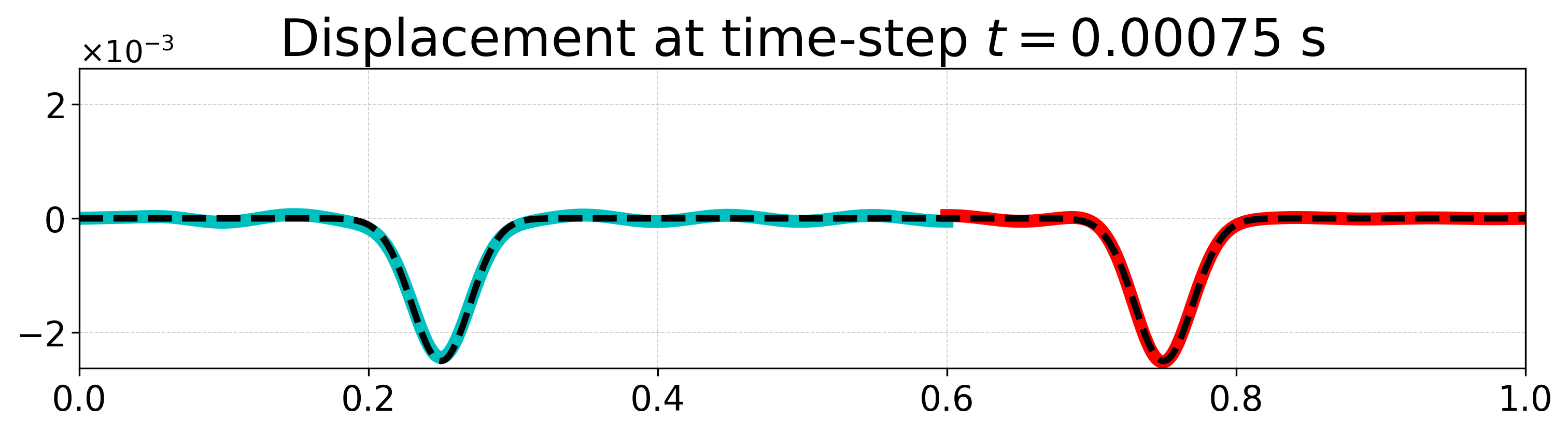}
    \end{subfigure}
    
    \vspace{1mm}
    
    \begin{subfigure}{0.45\linewidth}
        \centering
        \includegraphics[width=\linewidth]{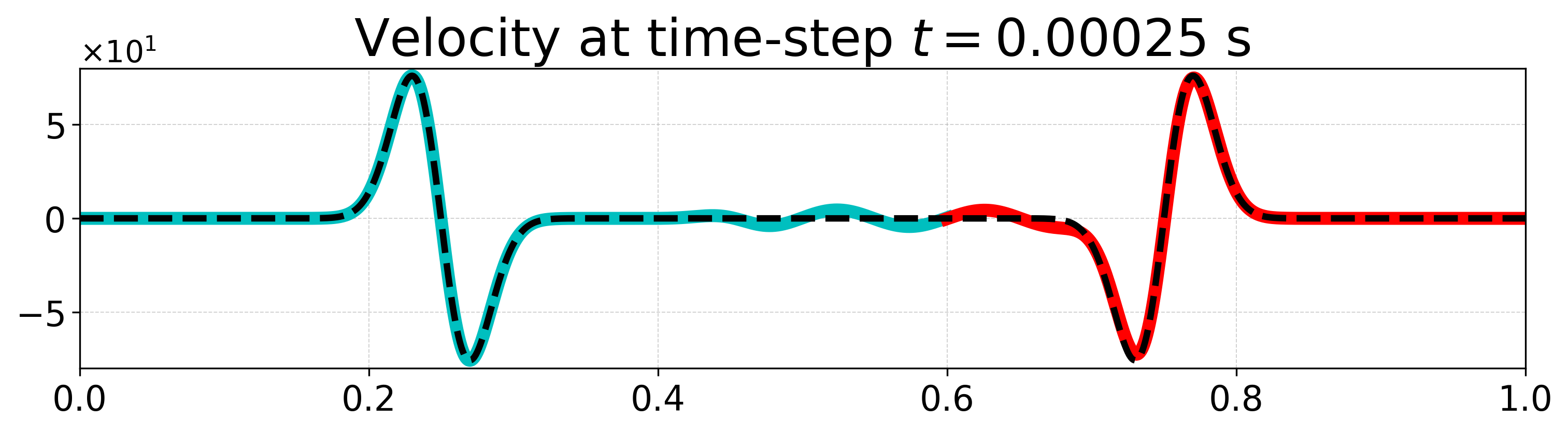}
    \end{subfigure}
    \hspace{3mm}
    \begin{subfigure}{0.45\linewidth}
        \centering
        \includegraphics[width=\linewidth]{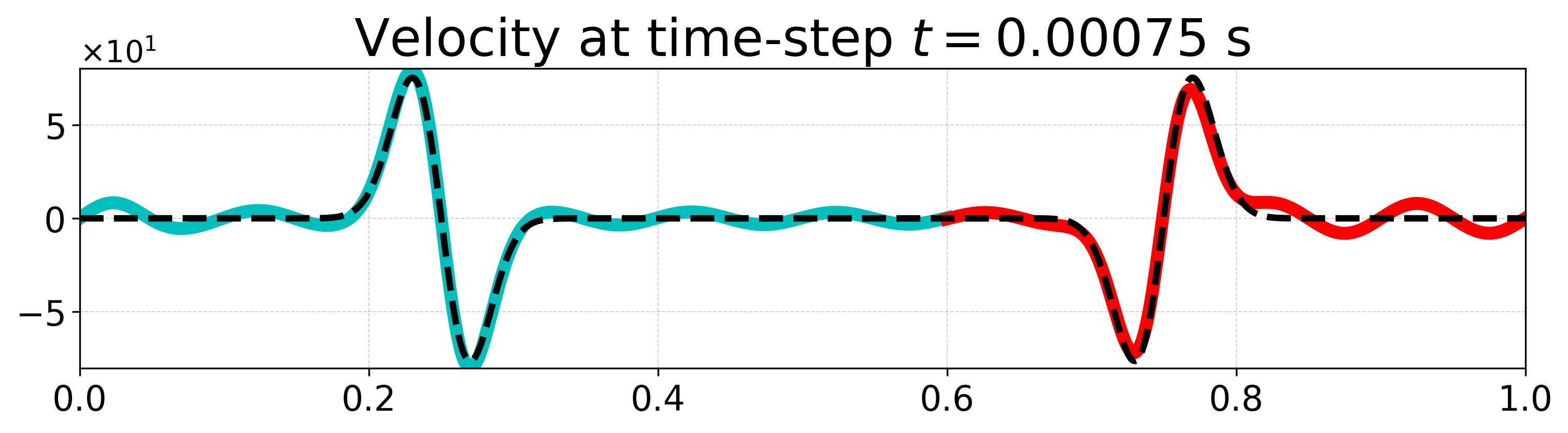}
    \end{subfigure}
    
    \vspace{1mm}
    
    \begin{subfigure}{0.45\linewidth}
        \centering
        \includegraphics[width=\linewidth]{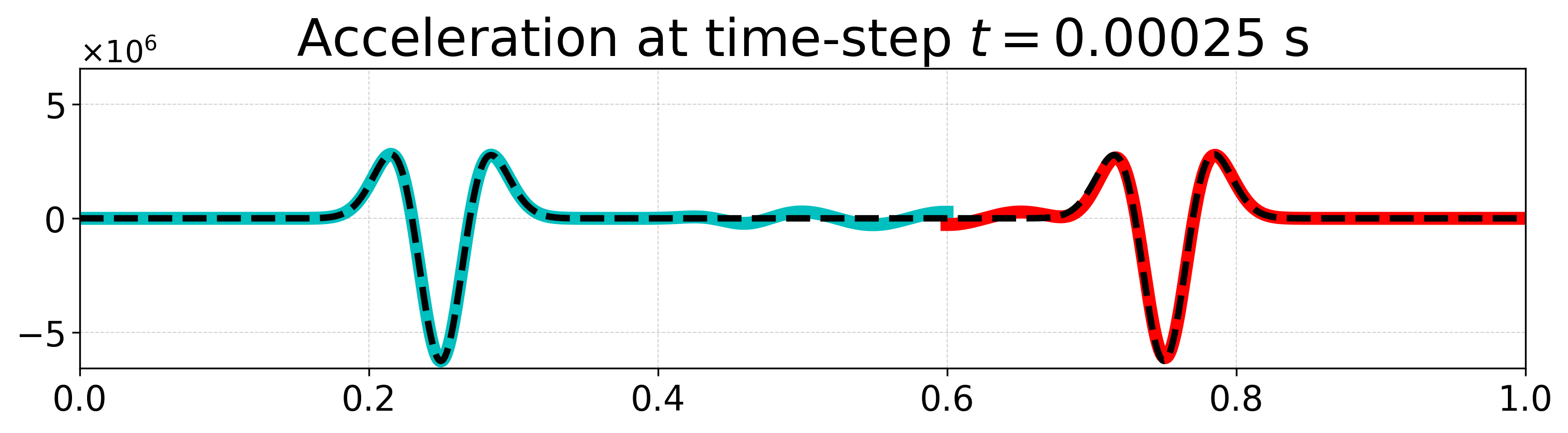}
        \caption{$t=2.5 \times 10^{-4} \, s$}
    \end{subfigure}
    \hspace{3mm}
    \begin{subfigure}{0.45\linewidth}
        \centering
        \includegraphics[width=\linewidth]{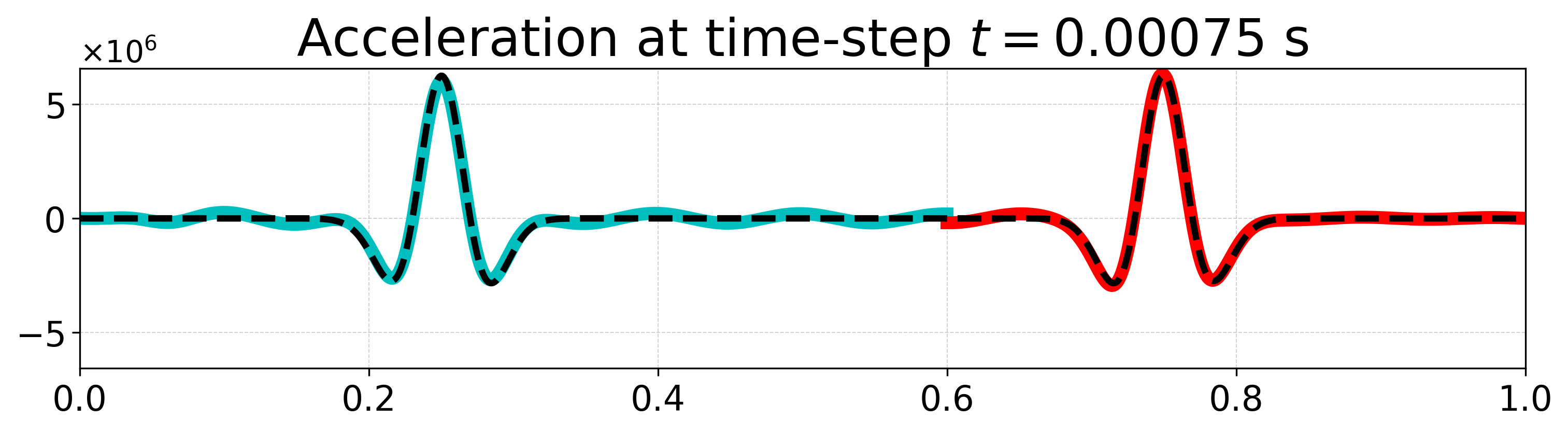}
        \caption{$t=7.5 \times 10^{-4} \, s$}
    \end{subfigure}
    \caption{FOM-FOM coupling: plots of the NO-SAM displacement, velocity, and acceleration solutions obtained with Robin–Robin transmission conditions (highest error parameter set, $\left[ \overline{\alpha}_{12}, \overline{\alpha}_{21}, \beta_{12},  \beta_{21} \right] = \left[ 10^{-1}, 1, 1, 3 \right]$) compared to the corresponding monolithic solution (dashed line).  The solutions computed for subdomains $\Omega_1$ and $\Omega_2$ are shown in blue and red, respectively.
}
    \label{fig:RR_Worst_plots}
\end{figure}

The displacement, velocity, and acceleration solutions at times $t=2.5 \times 10^{-4} \, s$ and $t=7.5 \times 10^{-4}\, s$ are shown for three key cases in Figures~\ref{fig:DN_plots} -- \ref{fig:RR_Worst_plots}. The solution obtained using the lowest error Robin parameters in Figure~\ref{fig:RR_Best_plots} is $25 \%$ more accurate than the alternating Dirichlet–Neumann result in Figure~\ref{fig:DN_plots}, while requiring $29\%$ fewer iterations. Both results exhibit a slight oscillation as the wave crosses the interface. This oscillation propagates backward along the bar, opposite to the direction of wave travel, and is most pronounced in the acceleration field. The result corresponding to the highest-error Robin parameters, shown in Figure~\ref{fig:RR_Worst_plots}, exhibits pronounced oscillations and discontinuities originating at the interface. Such oscillatory behavior would generally be unacceptable in engineering applications and highlights the importance of carefully tuning the $\alpha_{ij}$ and $\beta_{ij}$ parameters when employing Robin transmission conditions.

\section{NO-SAM for OpInf ROM coupling} \label{sec:schwarz_opinf}

Operator Inference \cite{willcox2016opinf} is a data-driven, non-intrusive, projection-based method for constructing low-dimensional surrogate models that approximate high-dimensional data. The data may be obtained from experiments or generated by numerical simulation; in this work, we focus on the latter. This section explores the use of NO-SAM to couple together subdomain-local OpInf ROMs with each other and with subdomain-local FOMs.  To the best of the authors' knowledge, this work is the first to explore the DD-based coupling of non-intrusive OpInf models using NO-SAM.

\subsection{Monolithic OpInf ROM overview} \label{sec:opinf}

Similar to intrusive projection-based ROM methodologies, OpInf begins by constructing a reduced basis from a given dataset using the proper orthogonal decomposition (POD) \cite{CJR:Holmes1996, CJR:Sirovich1987}. POD is a numerical technique that extracts a variance-maximizing set of features from a sequence of field snapshots. By capturing the most significant patterns in the data with a few dominant modes, POD reduces the complexity of the system, enabling the development of efficient, low-order surrogate models.

The dataset in this work comes from the monolithic numerical solution to the problem expressed in \eqref{eqn:matrix_vector_form}. Once we have snapshot data for $\tau>0$ separate states between $t_0$ and $t_f$, we can define the following snapshot matrix:
\begin{equation}
    \label{eqn:snapshot_matrix}
    ~U:=\left[ ~u(t_1), \, ~u(t_2), \, \dots \, , ~u(t_\tau)  \right] \in\mathbb{R}^{N \times\tau} 
\end{equation}

\noindent It is important to note that snapshot velocities $\dot{~u}(t_n)$ and accelerations $\ddot{~u}(t_n)$ can also be included in the snapshot matrix $~U$ along with the displacements, though this is not considered herein. The reduced basis is obtained by performing the singular value decomposition of $~U = ~\Phi ~\Sigma ~V^T$, where $~\Phi \in \mathbb{R}^{N \times N}$ is an orthogonal matrix whose first $R := \operatorname{rank}(~U) > 0$ columns span the column space of $~U$. The first $r \leq R$ columns of $~\Phi$ define an optimal $r$-dimensional basis for the columns of $~U$ in the $l^2$ sense. We denote this truncated basis by $~\Phi_r \in \mathbb{R}^{N \times r}$. The fraction of statistical energy captured by a basis of size $r$ is defined as
\begin{equation}
    \label{eqn:energy}
    E(r) \coloneq \frac{\sum_{i=1}^r \mu_i^2}{\sum_{i=1}^R \mu_i^2}
\end{equation}

\noindent where $\mu_i$ denotes the $i^{th}$ singular value of $~\Sigma$. 

For standard intrusive and projection-based ROMs, reduced operators are obtained by projecting the full order operators discretizing the governing PDEs onto the span of the reduced basis using $~\Phi$. The main disadvantage of this procedure is that it requires explicit access to the full order operators defining the governing problem. 
One selling point of OpInf ROMs is that the reduced operators are \emph{non-intrusively} inferred directly from the snapshot data used to build $~\Phi$, eliminating the need to access the full order system operators.

The first step in developing an appropriate OpInf ROM for our model problem (Section~\ref{sec:problem_setting}) is to establish the functional form of the FOM semi-discretization. As outlined in Section~\ref{sec:model}, applying the finite element method to discretize the weak form \eqref{eqn:weak_form} in space yields the functional representation given in \eqref{eqn:matrix_vector_form}. In the absence of body forces, the force vector reduces to


\begin{equation} \label{eq:Fi}
    F_i = \int_{\partial_T\Omega} ~T \cdot \boldsymbol{\xi}_i\ dS 
\end{equation}

\noindent for the $i^{\mathrm{th}}$ test function $\boldsymbol{\xi}_i$.  Our goal is to derive the functional form of a projection-based ROM corresponding to \eqref{eqn:matrix_vector_form}.   
Suppose the prescribed traction $~T:=~t$ is piecewise constant on the element faces along the subdomain interface $\Gamma$.  In this case, \eqref{eq:Fi} can be approximated as  $F_{i} =  ~t \cdot \int_{\Gamma} \boldsymbol{\xi}_i dS$, so that $~F$ takes the form: $~F = ~H ~t$ for a given prescribed traction vector, $~t$.   The resulting discrete system $~M \ddot{~u} + ~K ~u = ~H~t$ incorporates the traction, but not displacement boundary conditions.  Let $~i_d$ denote the solution 
indices on which the displacement
boundary conditions are imposed, and let $~i_u$ denote the remaining unconstrained indices.  
Suppose $~g$ denotes the vector of displacement values imposed at $~i_d$.  
Then, \eqref{eqn:matrix_vector_form} can 
be rewritten in the form 
\begin{equation} \label{eq:opinf_fom_with_bcs}
    \bar{~M}\, \ddot{\bar{~u}} + \bar{~K}\, \bar{~u} = \bar{~H}\, ~t + \bar{~B}\, ~g   
\end{equation}

\noindent where $\bar{~M}:= ~M \left[~i_u, ~i_u \right]$,   $\bar{~K}:= ~K \left[~i_u, ~i_u \right]$, $\bar{~u}: = ~u \left[ ~i_u \right]$, $\bar{~H}:= ~H \left[ ~i_u,: \right]$ and $\bar{~B}:=-~K \left[~i_u, ~i_d \right]$. 


Knowing that the structure of \eqref{eq:opinf_fom_with_bcs} is preserved for the intrusive projection-based ROM, we can assume a similar form for the semi-discretized OpInf formulation
\begin{equation}
    \label{eqn:opinf_matrix_vector}
    \hat{~M} \, \ddot{\hat{~u}} + \hat{~K} \, \hat{~u} = \hat{~H} \, ~t + \hat{~B} \, ~g
\end{equation}

\noindent \noindent where the ``hat" notation denotes reduced order operators. 
To learn the reduced order operators $\hat{~M}, \hat{~K} \in \mathbb{R}^{r \times r}$, $\hat{~H} \in \mathbb{R}^{r \times l}$, and $\hat{~B} \in \mathbb{R}^{r \times m}$ in the monolithic OpInf setting, we solve the following linear regression problem
\begin{equation}
    \label{eqn:monolithic_opinf_min}
    \min_{\hat{~M}, \hat{~K}, \hat{~H}, \hat{~B}} \sum_{p=1}^P \lVert \hat{~M} \, \ddot{\hat{~u}}(t_p) + \hat{~K} \, \hat{~u}(t_p) - \hat{~H} \, ~t(t_p) - \hat{~B} \, ~g(t_p) \rVert^2_2
\end{equation}

\noindent given boundary condition data $~t(t_p) \in \mathbb{R}^l$ and $~g(t_p) \in \mathbb{R}^m$. These terms are expressed as functions of time to account for cases where they are not constant. In such instances, snapshot data of the boundary conditions on the domain boundary $\partial\Omega$ must be collected at each time step. The reduced-order quantities $\hat{~u}(t_p) := ~\Phi_r^T ~u(t_p)$ and $\ddot{\hat{~u}}(t_p) := ~\Phi_r^T \ddot{~u}(t_p)$ are obtained from the displacement and acceleration snapshot matrices, respectively, where $P \leq \tau$ denotes the total number of snapshots taken from the FOM dataset. If acceleration snapshots are not available, they may be approximated using a finite-difference scheme. To avoid the trivial solution while solving \eqref{eqn:monolithic_opinf_min} in a fully unconstrained fashion—where each reduced operator identically equal to zero constitutes a valid minimum—we instead solve a modified linear regression problem obtained by applying the inverse of the reduced mass matrix to each term in \eqref{eqn:opinf_matrix_vector}
\begin{equation}
    \label{eqn:monolithic_opinf_min_inv_mass}
    \min_{\tilde{~K}, \tilde{~H}, \tilde{~B}} \sum_{p=1}^P \lVert \ddot{\hat{~u}}(t_p) + \tilde{~K} \, \hat{~u}(t_p) - \tilde{~H} \, ~t(t_p) - \tilde{~B} \, ~g(t_p) \rVert^2_2
\end{equation}

\noindent where the tilde notation denotes operators that have been premultiplied by the inverse of the reduced mass matrix, i.e.\ $\tilde{~K} := \hat{~M}^{-1} \hat{~K}$, $\tilde{~H} := \hat{~M}^{-1} \hat{~H}$, and $\tilde{~B} := \hat{~M}^{-1} \hat{~B}$. This formulation is employed in the numerical examples presented below. Once the reduced operators are determined, boundary conditions in the OpInf context are imposed directly by specifying $~g$ and $~t$ in \eqref{eqn:opinf_matrix_vector}. \\


\noindent {\bf Remark 4.2.  } While the assumption of constant traction at the interface is trivial in the 1D setting considered here, extending this formulation to higher dimensions may introduce additional challenges, and require alternate formulations/implementations.
One possible alternate formulation assumes that the quadrature leading to $~F$ in \eqref{eq:Fi} can be written as $~F = ~C~t'$ for some constant matrix $~C$ and a known (possibly nonlinear) function $~t'=~F(~t)$.  In this case, it is possible to learn reduced versions of $~C$ to construct the ROM.  In particular, we have 
\begin{equation} \label{eq:alt_rom}
   \hat{~M}\, \ddot{\hat{~u}} + \hat{~K}\,\hat{~u} = \hat{~C}\,~t' + \hat{~B}\,~g
\end{equation}
where $\hat{~C} := \boldsymbol{\Phi}_r^T~C$.  Learning the operators in \eqref{eq:alt_rom} is expected 
to lead to a viable ROM for each $~t$.  We plan to explore this research direction in future work.  
\\

An important consideration when implementing OpInf is the numerical stability of the resulting reduced order operators. 
One common technique to ensure the stability of the resulting reduced order operators is regularization, which can be applied to \eqref{eqn:monolithic_opinf_min_inv_mass} as
\begin{equation}
    \label{eqn:regularization}
    \min_{\tilde{~K}, \tilde{~H}, \tilde{~B}} 
    \sum_{p=1}^P \Big\lVert\ddot{ \hat{~u}}(t_p) + \tilde{~K} \, \hat{~u}(t_p) - \tilde{~H} \, ~t(t_p) - \tilde{~B} \, ~g(t_p) \Big\rVert_2^2 
    + \lambda^2 \Big( \lVert \tilde{~K} \rVert_F^2 + \lVert \tilde{~H} \rVert_F^2 + \lVert \tilde{~B} \rVert_F^2 \Big)
\end{equation}

\noindent where the regularization term serves to penalize large entries in the reduced order operators. This is particularly important for problems with a limited number of snapshots or when a large number of modes is required to capture the target energy of the system. A parameter sweep indicated that a regularization parameter of $\lambda = 10^{-4}$ provides an appropriate balance between stability and accuracy. Therefore, all numerical results presented herein use $\lambda = 10^{-4}$.

Another means of promoting the stability of the learned operators is to preprocess the dataset so that each variable has the same characteristic scale. This can involve variable transformations as well as scaling, shifting, and centering the data. It was found that scaling the traction data by a factor of $\frac{1}{\sigma_{\max}}$, where $\sigma_{\max}$ is the maximum stress defined in Section \ref{sec:schwarz_fom_results}, showed a significant improvement on the stability and accuracy of the resulting learned operators. This approach is applied in all relevant numerical examples going forward.

\subsection{NO-SAM formulation} \label{sec:schwarz_opinf_formulation}

As discussed in Section \ref{sec:intro}, this work focuses on the domain-decomposition-based coupling between full order and OpInf reduced order models. The generic setup consists of two non-overlapping subdomains, $\Omega_1$ and $\Omega_2$, with $\Omega_1 \, \cup \, \Omega_2 = \Omega$, as illustrated in Figure \ref{fig:non_overlapping_decomp}. In this context, we, for example, approximate the governing PDE in \eqref{eqn:weak_form} using a full order model in $\Omega_1$ and an OpInf reduced order model in $\Omega_2$. This formulation can naturally be extended to various combinations of both models, as well as decompositions that include more than two subdomains. 

A modified version of the OpInf minimization problem in \eqref{eqn:regularization} must be solved in the context of Schwarz coupling. In this case, each subdomain is associated with its own set of operators. Accordingly, the formulation in \eqref{eqn:regularization} is adapted for each subdomain $k$ as

\begin{equation}
    \label{eqn:regularization_schwarz}
    \min_{\tilde{~K}_k, \tilde{~H}_k, \tilde{~B}_k} 
    \sum_{p=1}^P \Big\lVert \ddot{\hat{~u}}_k(t_p) + \tilde{~K}_k \, \hat{~u}_k(t_p) - \tilde{~H}_k \, ~t_k(t_p) - \tilde{~B}_k \, ~g_k(t_p) \Big\rVert_2^2 
    + \lambda^2 \Big( \lVert \tilde{~K}_k \rVert_F^2 + \lVert \tilde{~H}_k \rVert_F^2 + \lVert \tilde{~B}_k \rVert_F^2 \Big)
\end{equation}

\noindent Where, for example, $\hat{~u}_k(t_p)$ denotes the reduced displacement within subdomain $\Omega_k$ at time $t_p$. For the case of Schwarz coupling, additional snapshot data containing the transmission conditions $~g$ and $~t$ on the interface $\Gamma$ at each time step is required. Only the operators required by the selected transmission condition for each subdomain are learned.

When Robin transmission conditions are employed, additional reduced order operators must be introduced to account for the extra stiffness and force contributions added in \eqref{eqn:robin_final}. The corresponding regularized OpInf minimization problem is formulated as
\begin{equation}
\label{eqn:opinf_robin}
\begin{aligned}
\min_{\substack{\tilde{~K}_k, \tilde{~S}_k, \tilde{~H}_k, \\ \tilde{~B}_k, \tilde{~R}_k}} \, 
& \sum_{p=1}^P \Big\lVert 
    \ddot{\hat{~u}}_k(t_p) 
    + \left[ \tilde{~K}_k + \tfrac{\beta_{ij}}{\alpha_{ij}} \tilde{~S}_k \right] \hat{~u}_k(t_p) 
    - \tilde{~H}_k \, ~t_k(t_p) 
    - \tilde{~B}_k \, ~g_k(t_p) 
    - \tfrac{1}{\alpha_{ij}}\tilde{~R}_k \, ~c_k(t_p) 
\Big\rVert_2^2 \\
& + \lambda^2 \Big( 
      \lVert \tilde{~K}_k \rVert_F^2 
    + \lVert \tilde{~S}_k \rVert_F^2 
    + \lVert \tilde{~H}_k \rVert_F^2 
    + \lVert \tilde{~B}_k \rVert_F^2 
    + \lVert \tilde{~R}_k \rVert_F^2 
\Big)
\end{aligned}
\end{equation}

\noindent where $\tilde{~S}_k$ denotes the reduced Robin stiffness operator, $\tilde{~R}_k$ is the reduced Robin boundary operator, and $~c_k(t_p):=\alpha _{ij} \, ~t_k(t_p) + \beta _{ij} \, ~g_k(t_p)$ represents the Robin transmission condition data prescribed on $\Gamma$ for subdomain $k$ at time $t_p$. The constants $\alpha_{ij}$ and $\beta_{ij}$ have been factored out of each operator, thereby removing their dependence on these parameters. By learning $\tilde{~S}_k$ and $\tilde{~K}_k$ independently, one can conveniently tune $\alpha_{ij}$ and $\beta_{ij}$ when necessary.
 For all relevant numerical examples, the Robin boundary condition data $~c_k(t_p)$ is scaled by a factor of $\frac{1}{\sigma_{\max}}$ to enhance stability.

We emphasize that in our OpInf NO-SAM algorithm, the reduced operators appearing in \eqref{eqn:opinf_robin} are learned once in the offline stage of the model reduction procedure, similar to the computation of the POD basis ${\bf \Phi}_r$.  Once these operators have been solved for,
 we can apply the online stage of the algorithm, in which time-stepping and subdomain synchronization via the NO-SAM algorithm is performed.  Specifically,
 in each time interval $I_n := [t_n, t_{n+1}]$,
 every Schwarz iteration involving the OpInf model is carried out according to the following steps:

\vspace{5pt} 
\begin{enumerate}
    \item Gather the transmission data from neighboring subdomains in the full order space. 
    \item Project the previous full order displacement, velocity, and acceleration snapshots into the reduced order basis using $~\Phi_r$.
    \item Advance the semi-discrete system in \eqref{eqn:opinf_matrix_vector} from $t_n$ to $t_{n+1}$ using the Newmark-$\beta$ scheme and the precomputed reduced order operators, to obtain the updated reduced displacements, velocities, and accelerations.
    \item Map the updated reduced solution back to the full order space.
    \item Impose Dirichlet transmission conditions strongly (if applicable) in the full order space.
    \item Check for convergence using \eqref{eqn:schwarz_tolerance}. If satisfied, advance to the next time step ($n \gets n+1$); otherwise, move to the next subdomain and return to step~1.
\end{enumerate}


\vspace{5pt} 

\subsection{Numerical results} \label{sec:schwarz_opinf_results}

We first evaluate the performance of various FOM–OpInf model combinations using both alternating Dirichlet–Neumann and Robin–Robin transmission conditions on a reproductive version of the 1D linear elastic wave propagation problem. For the Robin–Robin case, we employ the values of $\alpha_{ij}$ and $\beta_{ij}$ that correspond to the ``lowest error" values reported in Table \ref{tab:pareto_results}. When using the OpInf models, a total of 4000 snapshots of the corresponding monolithic displacement solution are collected over the interval $\left[ t_0, t_f \right]$, covering the entire reference monolithic simulation (meaning the results are reproductive rather than predictive). These snapshots are used to generate the POD basis. A total of $20$ and $17$ POD modes are required to capture $99.9\%$ of the snapshot energy for the left and right subdomains, respectively. Since the left subdomain is larger, it is expected that more modes are required to capture the same proportion of snapshot energy.

Table \ref{tab:opinf_fom_results} reports the error as defined in \eqref{eqn:error_avg}, the mean number of Schwarz iterations per time step required for convergence, and the CPU time for each FOM-OpInf model combination using both alternating Dirichlet-Neumann and Robin-Robin transmission conditions. The FOM-FOM results from Table \ref{tab:pareto_results} for both transmission condition types are included as a reference. For cases employing a subdomain-local OpInf model, the number of modes used to construct the corresponding POD basis is listed. Here, $M_1$ denotes the number of POD modes used for the left domain, and $M_2$ the number of POD modes used for the right domain. For POD bases retaining $99.9\%$ of the snapshot energy, it is clear that the inclusion of a subdomain OpInf model consistently increases both the number of Schwarz iterations required for convergence as well as the resulting error. To assess whether this behavior is due to insufficient modal content, we also consider a case with 34 and 29 POD modes for the left and right subdomains, respectively, corresponding to $99.999999\%$ energy retention in each subdomain.

\begin{table}[bt!]
    \centering
    \caption{FOM-OpInf and OpInf-OpInf coupling: comparison of Schwarz coupling strategies using Dirichlet–Neumann and Robin–Robin transmission conditions. Reported are the model types, reduced basis sizes $M_1$ and $M_2$, relative error, average iteration count, and CPU time.}
    \label{tab:opinf_fom_results}
    \begin{tabular}{lcccccc}
        \toprule[1pt]
        Transmission Type & Model & $M_1$ / $M_2$ & Error & Iterations & \makecell{CPU\\time $(s)$}\\
        \midrule[1pt]
        
        Dirichlet-Neumann & FOM-FOM & -- / -- & $3.43 \times 10^{-4}$ & $3.73$ & $52$ \\
        \addlinespace[2pt]\cmidrule(lr){1-6}\addlinespace[2pt]

        Dirichlet-Neumann & OpInf-FOM & 20 / -- & $7.37 \times 10^{-4}$ & $6.08$ & $64$ \\
        Dirichlet-Neumann & FOM-OpInf & -- / 17 & $7.31 \times 10^{-4}$ & $4.10$ & $47$ \\
        Dirichlet-Neumann & OpInf-OpInf & 20 / 17 & $9.63 \times 10^{-4}$ & $4.12$ & $37$ \\
        \addlinespace[2pt]\cmidrule(lr){1-6}\addlinespace[2pt]

        Dirichlet-Neumann & OpInf-FOM & 34 / -- & $1.29 \times 10^{-4}$ & $5.45$ & $60$ \\
        Dirichlet-Neumann & FOM-OpInf & -- / 29 & $1.53 \times 10^{-4}$ & $4.19$ & $51$ \\
        Dirichlet-Neumann & OpInf-OpInf & 34 / 29 & $2.03 \times 10^{-4}$ & $4.10$ & $39$ \\
        \addlinespace[2pt]\cmidrule(lr){1-6}\addlinespace[2pt]

        Robin-Robin &  FOM-FOM  & -- / -- & $2.57 \times 10^{-4}$ & $2.66$ & $39$ \\
        \addlinespace[2pt]\cmidrule(lr){1-6}\addlinespace[2pt]

        Robin-Robin &  OpInf-FOM  & 20 / -- & $6.35 \times 10^{-4}$ & $2.89$ & $32$ \\
        Robin-Robin &  FOM-OpInf  & -- / 17 & $2.18 \times 10^{-3}$ & $3.00$ & $35$ \\
        Robin-Robin &  OpInf-OpInf  & 20 / 17 & $2.47 \times 10^{-3}$ & $2.88$ & $26$ \\
        \addlinespace[2pt]\cmidrule(lr){1-6}\addlinespace[2pt]

        Robin-Robin &  OpInf-FOM  & 34 / -- & $1.62 \times 10^{-4}$ & $2.00$ & $27$ \\
        Robin-Robin &  FOM-OpInf  & -- / 29 & $9.42 \times 10^{23}$ & $3.80$ & $46$ \\
        Robin-Robin &  OpInf-OpInf & 34 / 29 & $1.22 \times 10^{-4}$ & $2.00$ & $22$ \\
        \bottomrule[1pt]
    \end{tabular}
\end{table}

\begin{figure}[b!]
    \centering
    
    \begin{subfigure}{0.45\linewidth}
        \centering
        \includegraphics[width=\linewidth]{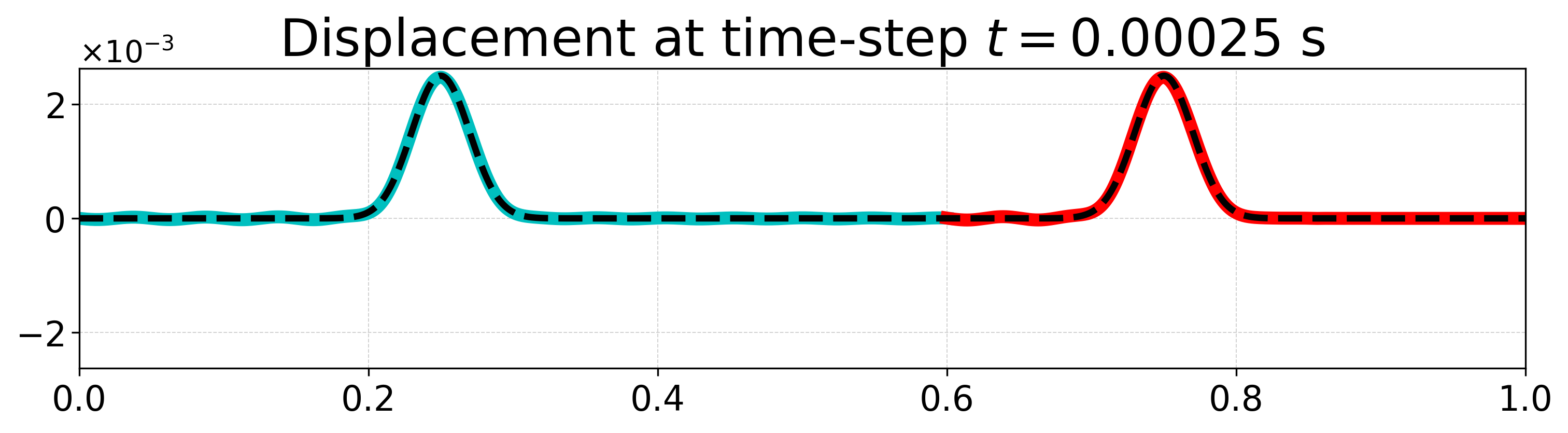}
    \end{subfigure}
    \hspace{3mm}
    \begin{subfigure}{0.45\linewidth}
        \centering
        \includegraphics[width=\linewidth]{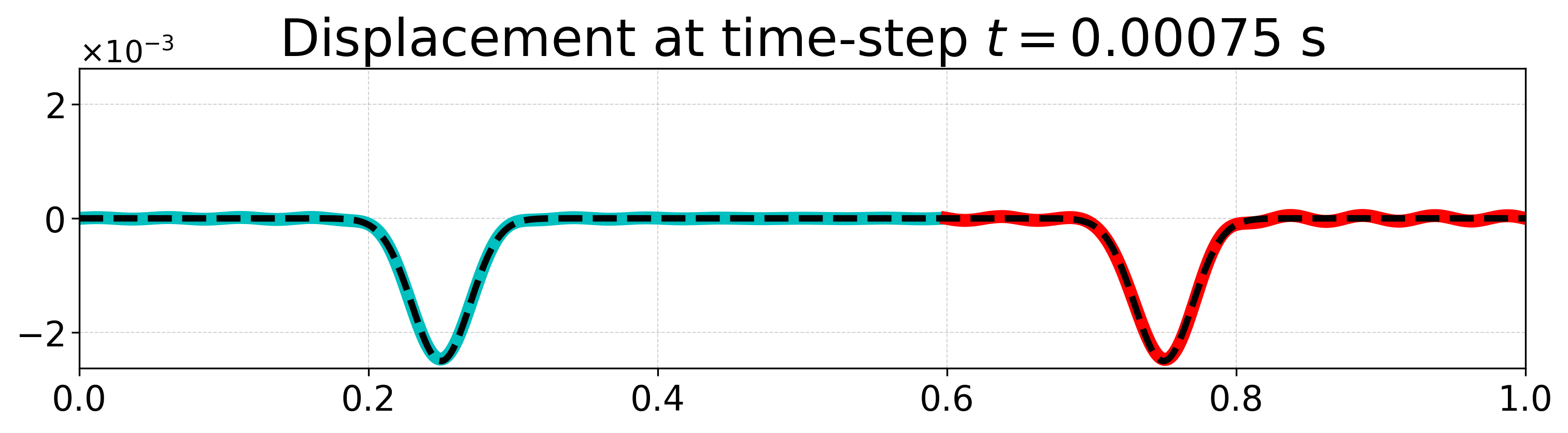}
    \end{subfigure}
    
    \vspace{1mm}
    
    \begin{subfigure}{0.45\linewidth}
        \centering
        \includegraphics[width=\linewidth]{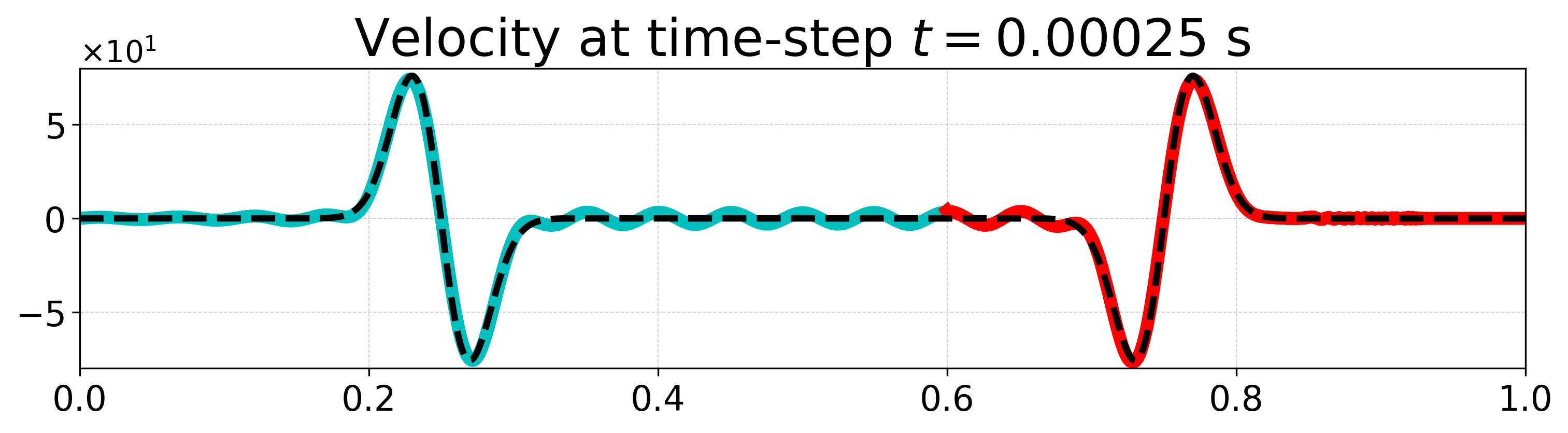}
    \end{subfigure}
    \hspace{3mm}
    \begin{subfigure}{0.45\linewidth}
        \centering
        \includegraphics[width=\linewidth]{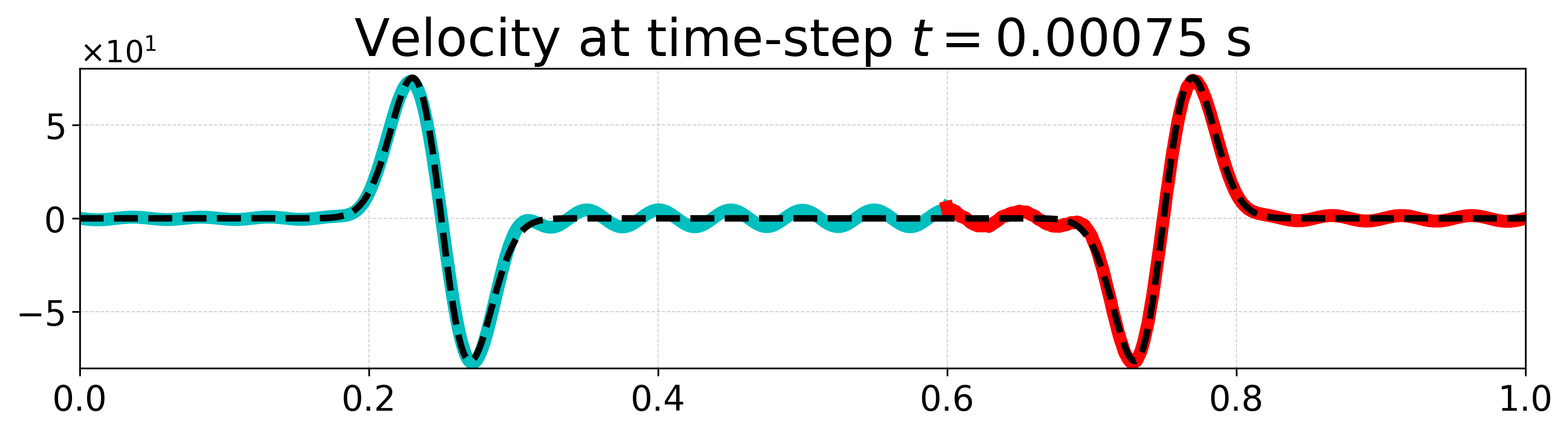}
    \end{subfigure}
    
    \vspace{1mm}
    
    \begin{subfigure}{0.45\linewidth}
        \centering
        \includegraphics[width=\linewidth]{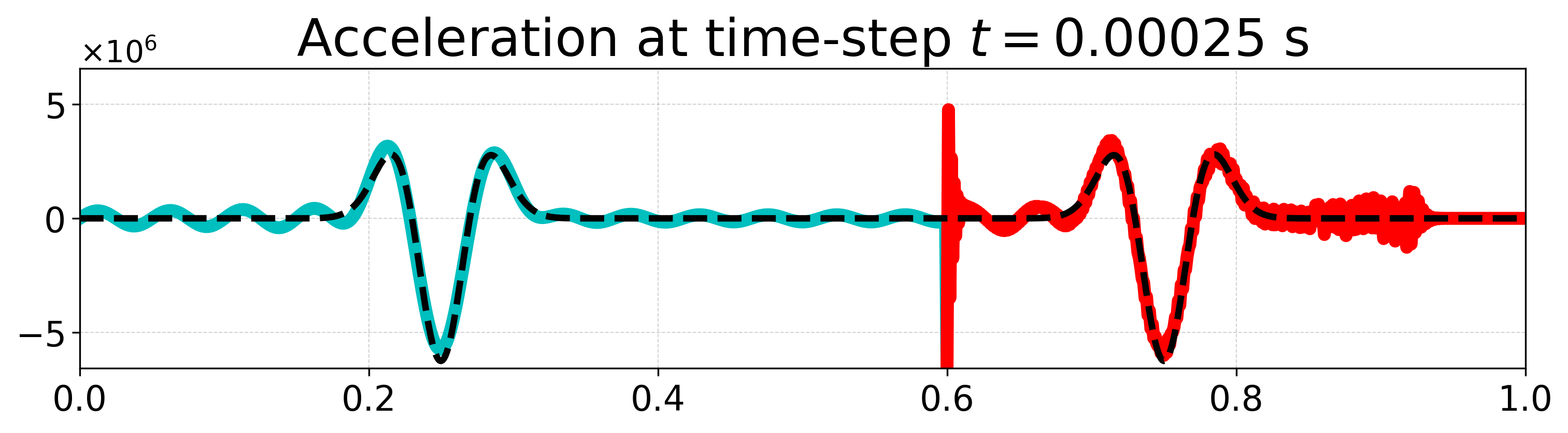}
        \caption{$t=2.5 \times 10^{-4} \, s$}
    \end{subfigure}
    \hspace{3mm}
    \begin{subfigure}{0.45\linewidth}
        \centering
        \includegraphics[width=\linewidth]{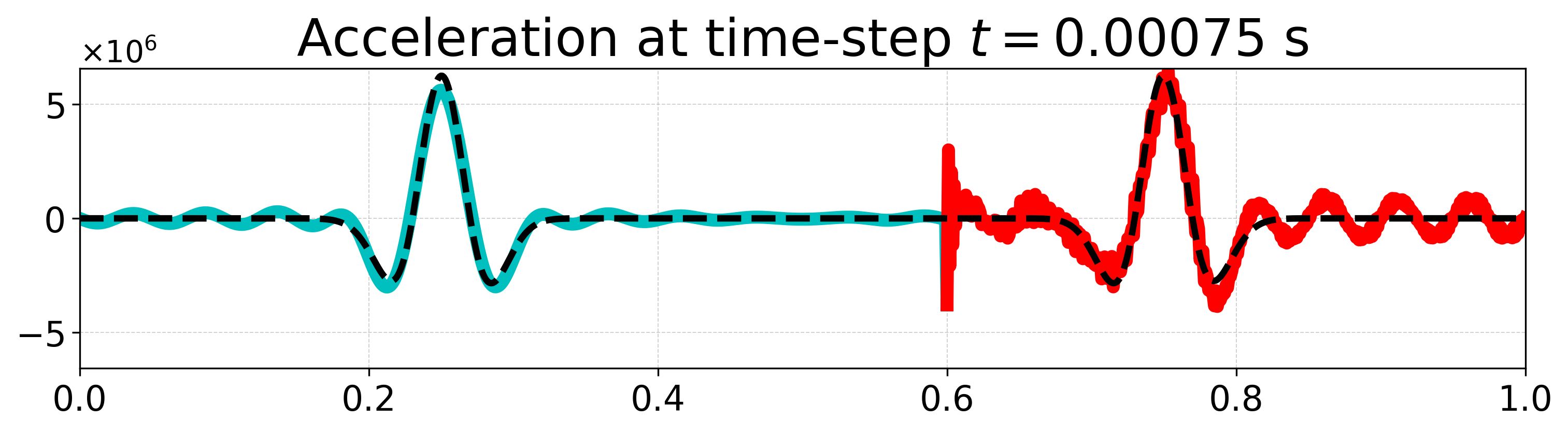}
        \caption{$t=7.5 \times 10^{-4} \, s$}
    \end{subfigure}
    \caption{OpInf-FOM coupling: plots of the NO-SAM displacement, velocity, and acceleration solutions obtained using alternating Dirichlet-Neumann transmission conditions, compared to the monolithic solution (dashed line). Here, the blue subdomain $\Omega_1$ employs an OpInf model and uses a Dirichlet transmission condition, and the red subdomain $\Omega_2$ employs a full order model and uses a Neumann transmission condition.}
    
    \label{fig:DN_OpInf_Instability}
\end{figure}

There are several interesting trends in Table \ref{tab:opinf_fom_results} that merit discussion. Focusing first on the alternating Dirichlet–Neumann results, the comparison between the lower and higher modal content settings indicates that the reduced accuracy in the former is primarily attributable to its limited modal content. However, for the OpInf-FOM configuration, an additional instability originating at the subdomain interface is present in both modal content settings, though it is significantly more severe in the lower modal content case. This instability, illustrated in Figure \ref{fig:DN_OpInf_Instability}, is most pronounced in the acceleration field, and likely contributes to the relatively large number of Schwarz iterations required for convergence and increase in CPU time in the OpInf-FOM case. More generally, the inclusion of a subdomain OpInf model can be expected to slow Schwarz convergence due to the accuracy loss inherent to ROMs. This increase is seen for both modal content settings. Despite this expected loss in accuracy, however, all three higher modal content settings actually exhibit a lower error than the reference FOM-FOM result. This can be attributed to the fact that, although OpInf yields a reduced order model, it employs basis functions that are optimal for representing the solution dataset, unlike the generic finite element basis. Consequently, the OpInf ROM can in some cases surpass the FOM in accuracy, as is also observed in \cite{SchwarzOpInfPaper:2025}.

\begin{figure}[t!]
    \centering
    
    \begin{subfigure}{0.45\linewidth}
        \centering
        \includegraphics[width=\linewidth]{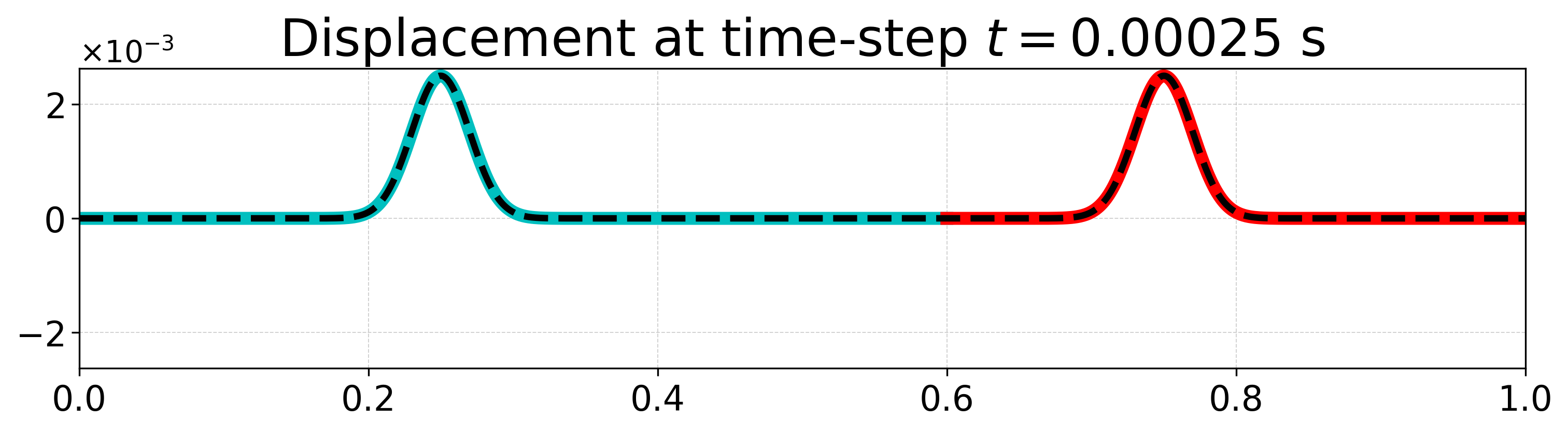}
    \end{subfigure}
    \hspace{3mm}
    \begin{subfigure}{0.45\linewidth}
        \centering
        \includegraphics[width=\linewidth]{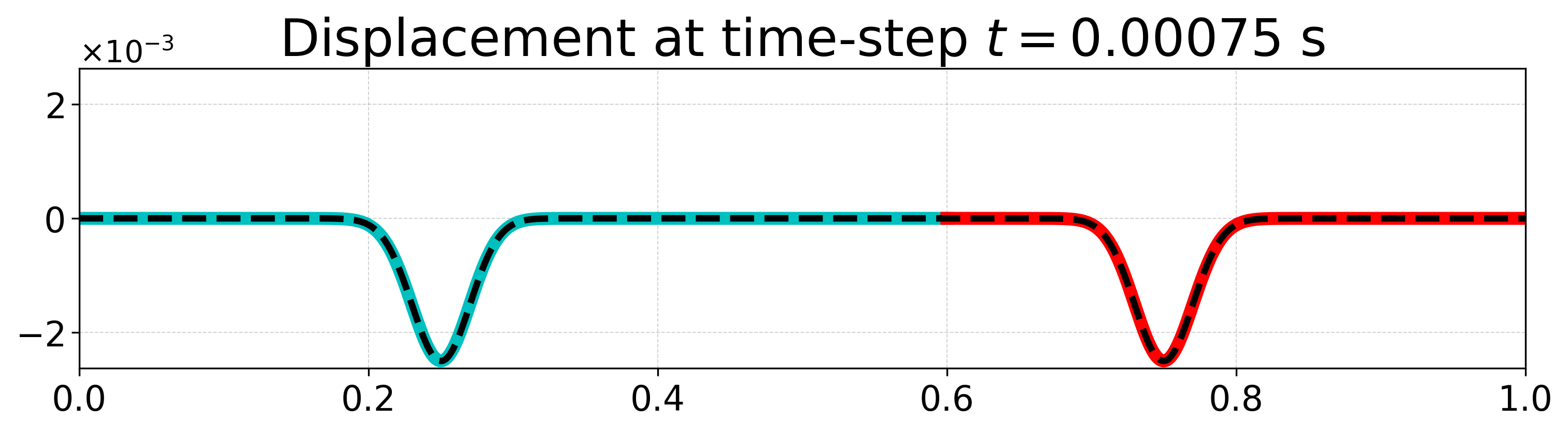}
    \end{subfigure}
    
    \vspace{1mm}
    
    \begin{subfigure}{0.45\linewidth}
        \centering
        \includegraphics[width=\linewidth]{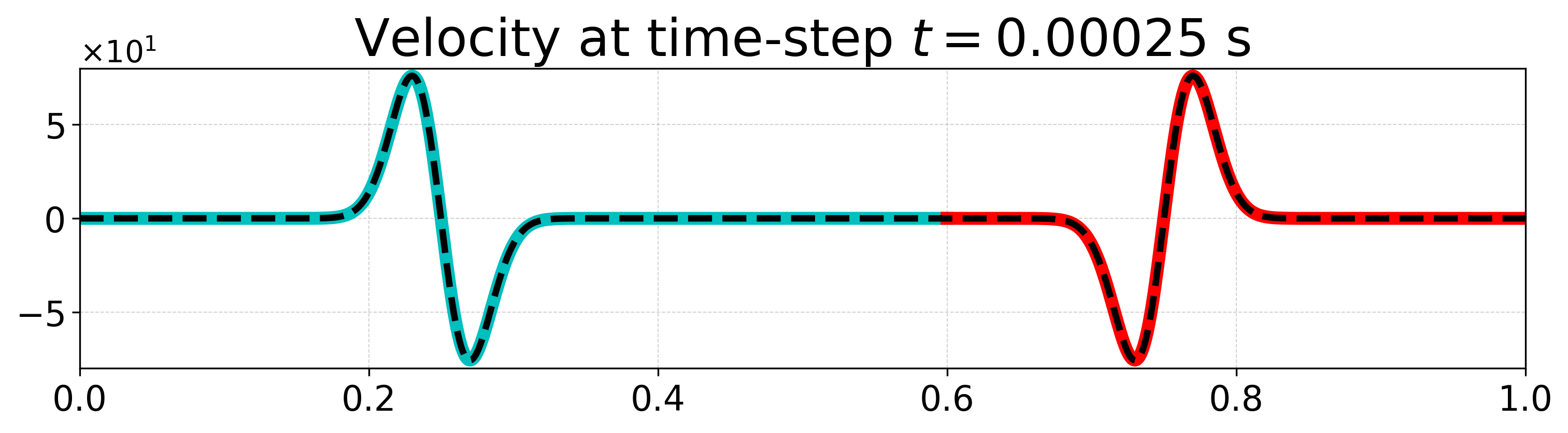}
    \end{subfigure}
    \hspace{3mm}
    \begin{subfigure}{0.45\linewidth}
        \centering
        \includegraphics[width=\linewidth]{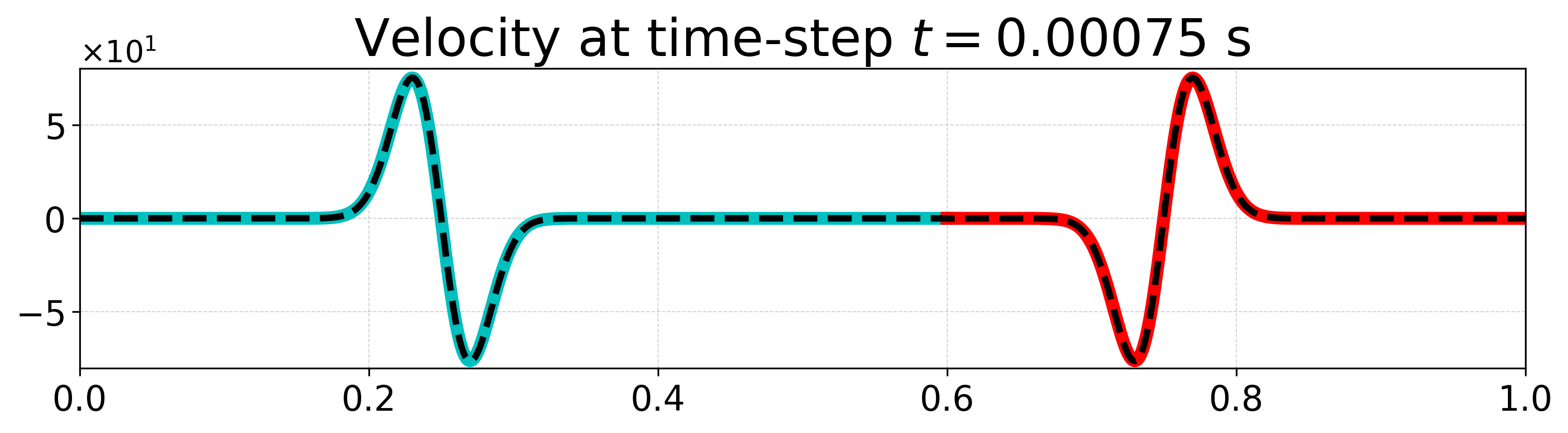}
    \end{subfigure}
    
    \vspace{1mm}
    
    \begin{subfigure}{0.45\linewidth}
        \centering
        \includegraphics[width=\linewidth]{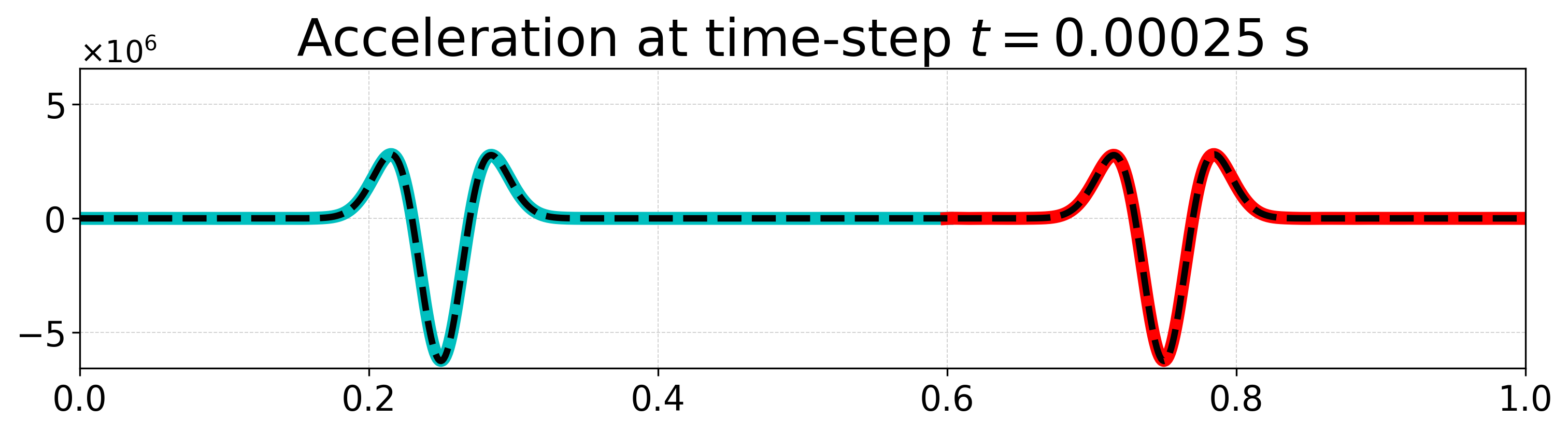}
        \caption{$t=2.5 \times 10^{-4} \, s$}
    \end{subfigure}
    \hspace{3mm}
    \begin{subfigure}{0.45\linewidth}
        \centering
        \includegraphics[width=\linewidth]{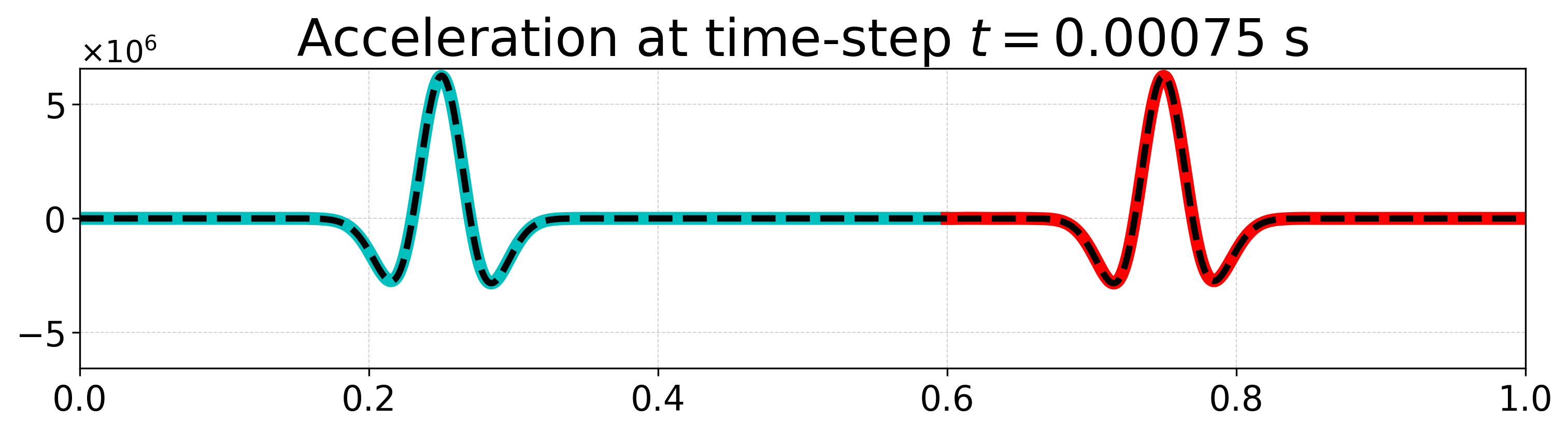}
        \caption{$t=7.5 \times 10^{-4} \, s$}
    \end{subfigure}
    \caption{OpInf-OpInf coupling: plots of the NO-SAM displacement, velocity, and acceleration solutions obtained using Robin-Robin transmission conditions (lowest error parameter set, $\left[ \overline{\alpha}_{12}, \overline{\alpha}_{21}, \beta_{12},  \beta_{21} \right] = \left[ 10^{-3}, 10^{-3}, 1, 1\right]$), compared to the monolithic solution (dashed line). Here, both the blue and red subdomains $\Omega_1$ and $\Omega_2$ employ an OpInf model.}
    
    \label{fig:RR_OpInf_Best}
\end{figure}

Similar trends are observed for the Robin–Robin case, where the large errors present in the lower modal content setting generally diminish as the modal content is increased. One exception is observed in the FOM-OpInf case, where the simulation becomes unstable. This instability is likely due to insufficient regularization, as the optimal regularization parameter $\delta$ may vary with the problem setting and model configuration. Among all tested cases, the OpInf-OpInf configuration with Robin–Robin transmission conditions and higher modal content yielded the best overall performance in terms of accuracy, iteration count, and CPU time. This simulation, illustrated in Figure \ref{fig:RR_OpInf_Best}, converged in two Schwarz iterations at every time step, achieving a speedup of $1.77\times$ and an accuracy improvement of $2.11\times$ compared to the reference FOM-FOM result. Recall that the $\alpha_{ij}$ and $\beta_{ij}$ parameters used for the Robin–Robin transmission condition were taken directly from the lowest-error parameter set identified in the FOM-FOM parameter sweep. It would, however, be interesting to generate a similar Pareto plot as in Figure~\ref{fig:FOM_FOM_robin_pareto} when coupling full order and reduced order models, to determine whether a similar trend emerges or if an even more optimal parameter set exists.

The final numerical example examines the accuracy and stability of the OpInf model within the predictive regime. The problem setting defined in Section \ref{sec:problem_setting} is used, with the simulation extended to $t_f=10^{-2} \, s$, corresponding to five periods of wave propagation along the bar. The corresponding monolithic displacement solution is used to collect 4000 snapshots over the original time interval, which spans half a period of wave propagation. The four configurations considered correspond to the combinations of alternating Dirichlet-Neumann and Robin-Robin transmission conditions with either full order or OpInf reduced order models. Here, both OpInf-OpInf configurations correspond to the higher modal content setting, with $M_1 = 34$ and $M_2 = 29$.

\begin{figure}
    \centering
    \includegraphics[width=1.0\linewidth]{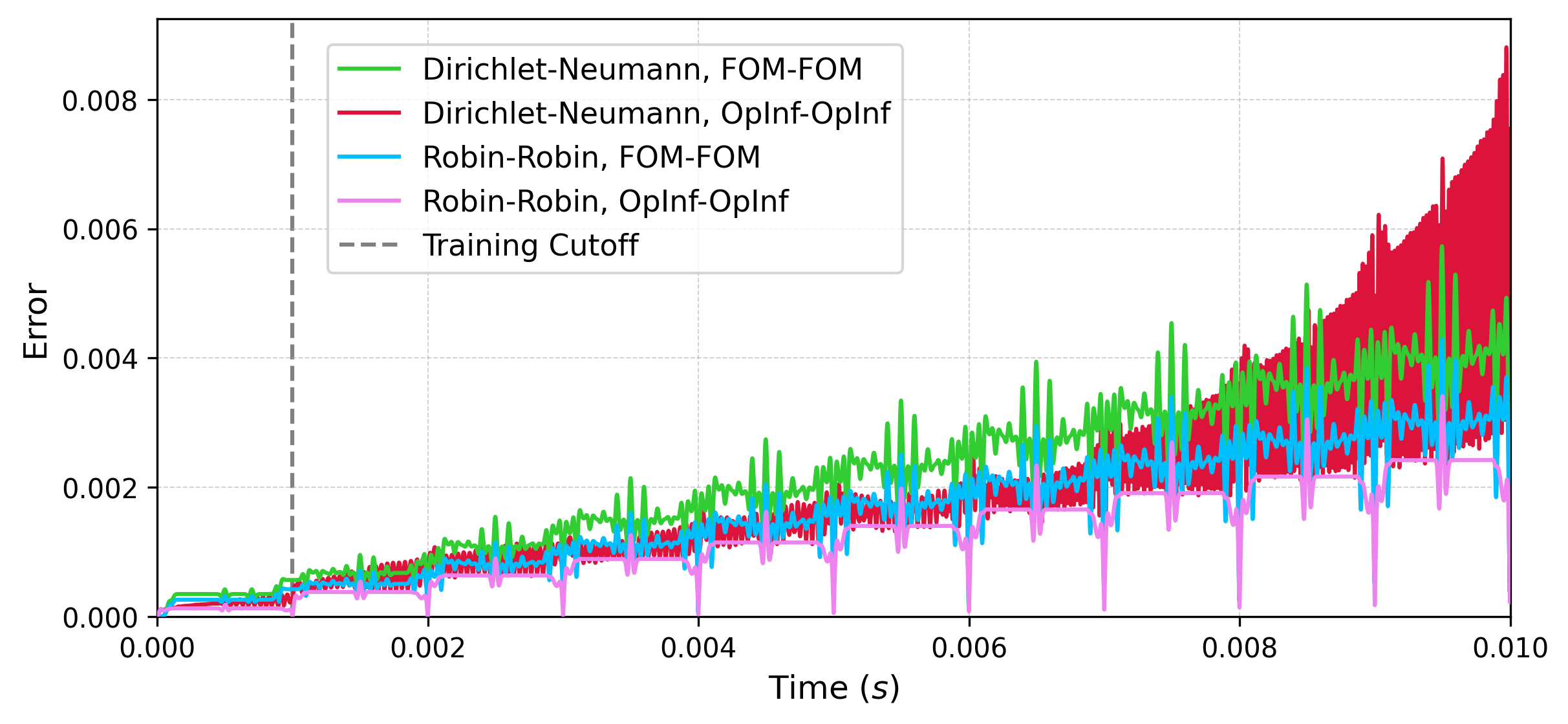}
   \caption{Error evolution over time for alternating Dirichlet-Neumann and Robin-Robin transmission conditions, using combinations of full order and OpInf reduced order models. The training cutoff, i.e., the final time step included in the training dataset, is marked by a vertical grey dashed line.}
    \label{fig:predictive_results}
\end{figure}

The total error, defined in \eqref{eqn:error_avg} without applying time-averaging, is plotted as a function of time for each case. The results are illustrated in Figure \ref{fig:predictive_results}. As expected, the Robin-Robin configurations generally exhibit lower error than their alternating Dirichlet-Neumann counterparts. Notably, the alternating Dirichlet-Neumann OpInf-OpInf configuration becomes unstable toward the end of the simulation, leading to a significant increase in error. This is likely related to the instability discussed previously for alternating Dirichlet-Neumann OpInf coupling, illustrated in Figure \ref{fig:DN_OpInf_Instability}. In contrast, both FOM-FOM configurations remain stable, though their errors grow steadily over time, likely due to spurious oscillations generated each time the wave crosses the interface, as discussed in Section~\ref{sec:schwarz_fom_results}. Among all cases, the Robin-Robin OpInf-OpInf configuration maintains the lowest error throughout the simulation and appears to be the least susceptible to instability.

\section{Conclusions and future work} \label{sec:conc}

This work has presented a first investigation into transmission conditions for NO-SAM applied to the coupling of subdomain-local FOMs and non-intrusive OpInf ROMs. Through systematic studies on a 1D linear elastic wave propagation problem, we demonstrated that the choice and tuning of interface conditions play a critical role in determining convergence, accuracy, and stability. Robin–Robin transmission conditions generally outperform alternating Dirichlet–Neumann formulations in terms of efficiency, provided that parameters are carefully selected. Our results also show that OpInf ROMs, when equipped with sufficient modal content, can deliver accuracy comparable to or even exceeding FOM–FOM couplings. At the same time, we observed interface instabilities and parameter sensitivities that highlight the importance of robust transmission design and regularization strategies.

Looking ahead, this study lays the groundwork for extending NO-SAM to more complex predictive regimes, multi-dimensional problems, and large-scale computational frameworks. In particular, integration into three-dimensional (3D) finite element codes such as {\tt Norma.jl} \cite{CJR:Norma} and deployment within production engineering codes such as {\tt SIERRA/SM} \cite{sierraSM56} will provide valuable opportunities for assessing scalability and practical impact. Moreover, the development of optimized domain decomposition strategies and hybrid FOM–ROM model libraries will be essential to fully harness the flexibility and efficiency of NO-SAM.
Finally, we plan to explore in future work the SAM-based coupling of alternate non-intrusive ROMs, e.g., ROMs constructed using reproducing kernel Hilbert spaces \cite{diaz2025interpretableflexiblenonintrusivereducedorder, diaz2025kernelmanifoldsnonlinearaugmentationdimensionality}.  
Together, these directions position the method as a promising tool for advancing scalable, non-intrusive coupling in computational mechanics and beyond.


\section*{Acknowledgments}

Support for this work was received through Sandia National Laboratories’ Laboratory Directed
Research and Development (LDRD) program and through the U.S. Department of Energy, Office of
Science, Office of Advanced Scientific Computing Research, Mathematical Multifaceted Integrated
Capability Centers (MMICCs) program, under Field Work Proposal 22025291 and the Multifaceted
Mathematics for Predictive Digital Twins (M2dt) project. 
C. Rodriguez acknowledges support from the U.S. Department of
Energy, Office of Science, Office of Advanced Scientific Computing Research, Department of
Energy Computational Science Graduate Fellowship under Award Number DE-SC0025528.
Additionally, the writing of this manuscript
was funded in part by I. Tezaur's Presidential Early Career Award for
Scientists and Engineers (PECASE).

Sandia National Laboratories is a multi-mission laboratory managed and operated by National Technology \& Engineering Solutions of Sandia, LLC (NTESS), a wholly owned subsidiary of Honeywell International Inc., for the U.S.~Department of Energy's National Nuclear Security Administration~(DOE/NNSA) under contract DE-NA0003525. This written work is authored by an employee of NTESS. The employee, not NTESS, owns the right, title and interest in and to the written work and is responsible for its contents. Any subjective views or opinions that might be expressed in the written work do not necessarily represent the views of the U.S.~Government.

\section*{Appendix}

Here, we derive the additional stiffness and boundary force operators arising from the Robin transmission conditions, as expressed in \eqref{eqn:robin_final}. Without loss of generality, consider the subdomain problem in $\Omega_1$. We begin with the force vector from the subdomain-local semi-discrete formulation, as given in the first line of \eqref{eqn:schwarz_1}. Here, we only consider the force contribution due to the transmission condition at the interface, denoted $~f^{(s+1)}_1$ and defined as

\begin{equation} \label{eq:force_appendix}
    ~f^{(s+1)}_1 := \int_{\Gamma} ~T^{(s+1)}_1 \cdot \boldsymbol{\xi}\ dS 
\end{equation}

\noindent It is important to note that this term is only integrated over the subdomain interface $\Gamma$. Solving for $~T^{(s+1)}_1$ in the third line of \eqref{eqn:schwarz_1} and substituting the result into \eqref{eq:force_appendix}, while using the expression for $~\lambda^{(s+1)}_1$ and assuming no relaxation ($\theta_1 = 1$), yields

\begin{equation} \label{eq:force_appendix_sub_robin}
    ~f^{(s+1)}_1 = \int_{\Gamma} \left[ \frac{1}{\alpha_{12}} ~c_1^{(s)} - \frac{\beta_{12}}{\alpha_{12}}~u_1^{(s+1)}\right] \cdot \boldsymbol{\xi} \ dS 
\end{equation}

\noindent where $~c_1^{(s)} := \alpha_{12} ~T_2^{(s)} + \beta_{12} ~u_2^{(s)}$ denotes the Robin transmission condition prescribed by the neighboring subdomain. The force term in \eqref{eq:force_appendix_sub_robin} can then be decomposed into two separate integrals:

\begin{equation}
    \label{eq:define_S}
        ~S^{(s+1)}_1 := \frac{\beta_{12}}{\alpha_{12}}\int_{\Gamma} ~u_1^{(s+1)} \cdot \boldsymbol{\xi} \, dS
\end{equation}

\begin{equation}
    \label{eq:define_R}
        ~R^{(s+1)}_1 := \frac{1}{\alpha_{12}}\int_{\Gamma}  ~c_1^{(s)} \cdot \boldsymbol{\xi} \, dS
\end{equation}

\noindent Together, these expressions, when substituted into the external force vector of the first line of \eqref{eqn:schwarz_1}, yield the final semi-discrete matrix-vector form. Here, $~S^{(s+1)}_1$ contributes an additional stiffness component, while $~R^{(s+1)}_1$ contributes an additional external force component. The remaining $~F^{(s+1)}_1$ vector in \eqref{eqn:robin_final} represents the external body force together with any applied traction on the outer domain boundary $\partial_{T} \Omega$. Observe that the negative sign in \eqref{eq:force_appendix_sub_robin} disappears when the second term is moved to the left-hand side of \eqref{eqn:robin_final}.


\bibliographystyle{siam}
\bibliography{CameronRodriguez}

\begin{thebibliography}{10}

\bibitem{Barnett:2022Schwarz}
{\sc J.~Barnett, I.~Tezaur, and A.~Mota}, {\em The {S}chwarz alternating method
  for the seamless coupling of nonlinear reduced order models and full order
  models}.
\newblock ArXiv pre-print, 2022.

\bibitem{sierraSM56}
{\sc G.~L. Bergel, F.~B. Beckwith, G.~J. de~Frias, K.~M. Manktelow, M.~T.
  Merewether, S.~T. Miller, K.~J. Parmar, T.~S. Shelton, J.~T. Thomas, J.~T.
  Trageser, B.~T. Treweek, M.~V. Veilleux, and E.~B. Wagman}, {\em
  {Sierra/SolidMechanics 5.6 User's Manual}}, Tech. Rep. SAND2022-3705, Sandia
  National Laboratories, Mar. 2022.

\bibitem{Bochev:2024}
{\sc P.~Bochev, J.~Owen, P.~Kuberry, and J.~Connors}, {\em Dynamic flux
  surrogate-based partitioned methods for interface problems}, Computer Methods
  in Applied Mechanics and Engineering, 429 (2024), p.~117115.

\bibitem{Chung:2024}
{\sc S.~W. Chung, Y.~Choi, P.~Roy, T.~Moore, T.~Roy, T.~Y. Lin, D.~T. Nguyen,
  C.~Hahn, E.~B. Duoss, and S.~E. Baker}, {\em Train small, model big: Scalable
  physics simulators via reduced order modeling and domain decomposition},
  Computer Methods in Applied Mechanics and Engineering, 427 (2024), p.~117041.

\bibitem{Cinquegrana:2011}
{\sc D.~Cinquegrana, A.~Viviani, and R.~Donelli}, {\em {A hybrid method based
  on POD and domain decomposition to compute the 2-D aerodynamic flow field -
  incompressible validation}}, 2011, pp.~AIMETA 2011-- XX Congresso
  dell'Associazione Italiana di Meccanica Teorica e Applicata, Bologna, ITA.

\bibitem{Corigliano:2015}
{\sc A.~Corigliano, M.~Dossi, and S.~Mariani}, {\em Model order reduction and
  domain decomposition strategies for the solution of the dynamic
  elastic–plastic structural problem}, Computer Methods in Applied Mechanics
  and Engineering, 290 (2015), pp.~127--155.

\bibitem{deCastro:2023}
{\sc A.~{de Castro}, P.~Bochev, P.~Kuberry, and I.~Tezaur}, {\em Explicit
  synchronous partitioned scheme for coupled reduced order models based on
  composite reduced bases}, Computer Methods in Applied Mechanics and
  Engineering,  (2023), p.~116398.

\bibitem{deCastro:2022}
{\sc A.~de~Castro, P.~Kuberry, I.~Tezaur, and P.~Bochev}, {\em A novel
  partitioned approach for reduced order model -- finite element model
  ({ROM}-{FEM}) and {ROM}-{ROM} coupling}.
\newblock ArXiv pre-print.

\bibitem{diaz2025interpretableflexiblenonintrusivereducedorder}
{\sc A.~N. Diaz, S.~A. McQuarrie, J.~T. Tencer, and P.~J. Blonigan}, {\em
  Interpretable and flexible non-intrusive reduced-order models using
  reproducing kernel hilbert spaces}, 2025.

\bibitem{diaz2025kernelmanifoldsnonlinearaugmentationdimensionality}
{\sc A.~N. Diaz, J.~T. Needels, I.~K. Tezaur, and P.~J. Blonigan}, {\em Kernel
  manifolds: nonlinear-augmentation dimensionality reduction using reproducing
  kernel hilbert spaces}, 2025.

\bibitem{Discacciati:2024}
{\sc N.~Discacciati and J.~S. Hesthaven}, {\em Model reduction of coupled
  systems based on non-intrusive approximations of the boundary response maps},
  Computer Methods in Applied Mechanics and Engineering, 420 (2024), p.~116770.

\bibitem{Ebrahimi:2024}
{\sc M.~Ebrahimi and M.~Yano}, {\em A hyperreduced reduced basis element method
  for reduced-order modeling of component-based nonlinear systems}, Computer
  Methods in Applied Mechanics and Engineering, 431 (2024), p.~117254.

\bibitem{Farcas:2023}
{\sc I.-G. Farcas, R.~P. Gundevia, R.~Munipalli, and K.~E. Willcox}, {\em
  Improving the accuracy and scalability of large-scale physics-based
  data-driven reduced modeling via domain decomposition}.
\newblock ArXiv pre-print, 2023.

\bibitem{SandiaLabNews}
{\sc M.~Fisher}, {\em {Automating complex 3D modeling}}.
\newblock Sandia Lab News, 2020.

\bibitem{Gander:2006}
{\sc M.~J. Gander}, {\em Optimized schwarz methods}, SIAM Journal on Numerical
  Analysis, 44 (2006), pp.~699--731.

\bibitem{Gkimisis:2025}
{\sc L.~Gkimisis, N.~Aretz, M.~Tezzele, T.~Richter, P.~Benner, and K.~E.
  Willcox}, {\em Non-intrusive reduced-order modeling for dynamical systems
  with spatially localized features}, 2025.

\bibitem{Hawkins:2024}
{\sc E.~Hawkins, P.~Kuberry, and P.~Bochev}, {\em An optimization-based
  coupling of reduced order models with efficient reduced adjoint basis
  generation approach}.
\newblock ArXiv pre-print, 2024.

\bibitem{Hoang:2021}
{\sc C.~Hoang, Y.~Choi, and K.~Carlberg}, {\em Domain-decomposition
  least-squares petrov–galerkin (dd-lspg) nonlinear model reduction},
  Computer Methods in Applied Mechanics and Engineering, 384 (2021), p.~113997.

\bibitem{CJR:Holmes1996}
{\sc P.~Holmes, J.~Lumley, and G.~Berkooz}, {\em Turbulence, Coherent
  Structures, Dynamical Systems and Symmetry}, Cambridge University Press,
  1996.

\bibitem{Huang:2024}
{\sc C.~Huang and K.~Duraisamy}, {\em Predictive reduced order modeling of
  chaotic multi-scale problems using adaptively sampled projections}, Journal
  of Computational Physics, 491 (2023), p.~112356.

\bibitem{hughes1987finite}
{\sc T.~J.~R. Hughes}, {\em The Finite Element Method: Linear Static and
  Dynamic Finite Element Analysis}, Prentice-Hall, Englewood Cliffs, NJ, 1987.

\bibitem{Iollo:2022}
{\sc A.~Iollo, G.~Sambataro, and T.~Taddei}, {\em A one-shot overlapping
  {S}chwarz method for component-based model reduction: application to
  nonlinear elasticity}, Computer Methods in Applied Mechanics and Engineering,
  404 (2023), p.~115786.

\bibitem{Kerfriden:2013}
{\sc P.~Kerfriden, O.~Goury, T.~Rabczuk, and S.~Bordas}, {\em A partitioned
  model order reduction approach to rationalise computational expenses in
  nonlinear fracture mechanics}, Computer Methods in Applied Mechanics and
  Engineering, 256 (2013), pp.~169--188.

\bibitem{Kerfriden:2012}
{\sc P.~Kerfriden, J.~C. Passieux, and S.~P.~A. Bordas}, {\em Local/global
  model order reduction strategy for the simulation of quasi-brittle fracture},
  International Journal for Numerical Methods in Engineering, 89 (2012),
  pp.~154--179.

\bibitem{CJR:Norma}
{\sc S.~N. Laboratories}, {\em {Norma.jl}: A julia testbed for coupling and
  multiphysics}.
\newblock \url{https://github.com/sandialabs/Norma.jl}, 2025.
\newblock Commit on GitHub: accessed on September 2, 2025.

\bibitem{LiD3M}
{\sc K.~Li, K.~Tang, T.~Wu, and Q.~Liao}, {\em D3m: A deep domain decomposition
  method for partial differential equations}, IEEE Access, 8 (2020),
  pp.~5283--5294.

\bibitem{LiDeepDDM}
{\sc W.~Li, X.~Xiang, and Y.~Xu}, {\em Deep domain decomposition method:
  {E}lliptic problems}, Proceedings of Machine Learning Research, 107 (2020),
  pp.~269--286.

\bibitem{Lions:1990}
{\sc P.~Lions}, {\em {On the Schwarz alternating method III: A variant for
  nonoverlapping subdomains}}.
\newblock in Third International Symposium on Domain Decomposition Methods for
  Partial Differential Equations, T. F. Chan, R. Glowinski, J. Periaux, and O.
  B. Widlund, eds., Society for Industrial and Applied Mathematics, 1990.

\bibitem{Maier:2014}
{\sc I.~Maier and B.~Haasdonk}, {\em {A {D}irichlet-{N}eumann reduced basis
  method for homogeneous domain decomposition problems}}, Applied Numerical
  Mathematics, 78 (2014), pp.~31--48.

\bibitem{Moore:2024}
{\sc I.~Moore, C.~Wentland, A.~Gruber, and I.~Tezaur}, {\em Domain
  decomposition-based coupling of operator inference reduced order models via
  the schwarz alternating method}.
\newblock ArXiV pre-print, 2024.

\bibitem{Mota:2025}
{\sc A.~Mota, D.~Koliesnikova, I.~Tezaur, and J.~Hoy}, {\em A fundamentally new
  coupled approach to contact mechanics via the dirichlet-neumann schwarz
  alternating method}, International Journal for Numerical Methods in
  Engineering, 126 (2025), p.~e70039.

\bibitem{Mota:2017}
{\sc A.~Mota, I.~Tezaur, and C.~Alleman}, {\em The {S}chwarz alternating method
  in solid mechanics}, Computer Methods in Applied Mechanics and Engineering,
  319 (2017), pp.~19--51.

\bibitem{Mota:2022}
{\sc A.~Mota, I.~Tezaur, and G.~Phlipot}, {\em {The {S}chwarz alternating
  method for dynamic solid mechanics}}, Int. J. Numer. Meth. Engng,  (2022),
  pp.~1--36.

\bibitem{SchwarzOpInfPaper:2025}
{\sc E.~Parish, A.~Gruber, I.~Tezaur, I.~Moore, C.~Wentland, and A.~Mota}, {\em
  {Subdomain-Local Coupling with Operator Inference and the Schwarz Alternating
  Method}}.
\newblock in preparation, 2025.

\bibitem{willcox2016opinf}
{\sc B.~Peherstorfer and K.~Willcox}, {\em Data-driven operator inference for
  nonintrusive projection-based model reduction}, Computer Methods in Applied
  Mechanics and Engineering, 306 (2016), pp.~196--215.

\bibitem{Prusak:2023}
{\sc I.~Prusak, M.~Nonino, D.~Torlo, F.~Ballarin, and G.~Rozza}, {\em An
  optimisation–based domain–decomposition reduced order model for the
  incompressible {N}avier-{S}tokes equations}, Computers \& Mathematics with
  Applications, 151 (2023), pp.~172--189.

\bibitem{Prusak:2024}
{\sc I.~Prusak, D.~Torlo, M.~Nonino, and G.~Rozza}, {\em An
  optimisation–based domain–decomposition reduced order model for
  parameter–dependent non–stationary fluid dynamics problems}, Computers \&
  Mathematics with Applications, 166 (2024), pp.~253--268.

\bibitem{Radermacher:2014}
{\sc A.~Radermacher and S.~Reese}, {\em Model reduction in elastoplasticity:
  proper orthogonal decomposition combined with adaptive sub-structuring},
  Computational Mechanics, 54 (2014), pp.~677--687.

\bibitem{Schwarz:1870}
{\sc H.~A. Schwarz}, {\em Ueber einen Grenz{\"u}bergang durch alternirendes
  Verfahren}, Z{\"u}rcher u. Furrer, 1870.

\bibitem{CJR:Sirovich1987}
{\sc L.~Sirovich}, {\em Turbulence and the dynamics of coherent structures,
  part iii: dynamics and scaling}, Q. Appl. Math., 45 (1987), pp.~583--590.

\bibitem{Smetana:2023}
{\sc K.~Smetana and T.~Taddei}, {\em Localized model reduction for nonlinear
  elliptic partial differential equations: Localized training, partition of
  unity, and adaptive enrichment}, SIAM Journal on Scientific Computing, 45
  (2023), pp.~A1300--A1331.

\bibitem{Snyder:2023}
{\sc W.~Snyder, I.~Tezaur, and C.~Wentland}, {\em Domain decomposition-based
  coupling of physics-informed neural networks via the {S}chwarz alternating
  method}.
\newblock ArXiv pre-print, 2023.

\bibitem{Taddei:2024}
{\sc T.~Taddei, X.~Xu, and L.~Zhang}, {\em A non-overlapping optimization-based
  domain decomposition approach to component-based model reduction of
  incompressible flows}, Journal of Computational Physics, 509 (2024),
  p.~113038.

\bibitem{Wang:2022}
{\sc H.~Wang, R.~Planas, A.~Chandramowlishwaran, and R.~Bostanabad}, {\em
  Mosaic flows: A transferable deep learning framework for solving pdes on
  unseen domains}, Computer Methods in Applied Mechanics and Engineering, 389
  (2022), p.~114424.

\bibitem{Wang:2025}
{\sc W.~Wang, M.~Hakimzadeh, H.~Ruan, and S.~Goswami}, {\em Accelerating
  multiscale modeling with hybrid solvers: Coupling fem and neural operators
  with domain decomposition}, 2025.

\bibitem{Wentland:2025}
{\sc C.~R. Wentland, F.~Rizzi, J.~L. Barnett, and I.~K. Tezaur}, {\em The role
  of interface boundary conditions and sampling strategies for {S}chwarz-based
  coupling of projection-based reduced order models}, Journal of Computational
  and Applied Mathematics, 465 (2025), p.~116584.

\bibitem{Zanolli:1987}
{\sc P.~Zanolli}, {\em Domain decomposition algorithms for spectral methods},
  CALCOLO, 24 (1987), pp.~201--240.

\end{thebibliography}

\end{document}